\documentclass[11pt]{article}

\usepackage{amsmath}
\usepackage{amssymb}
\usepackage{amsthm}
\usepackage[margin=1in]{geometry}
\usepackage{url}
\usepackage[colorlinks = true,
            linkcolor = blue,
            urlcolor  = blue,
            citecolor = blue,
            anchorcolor = blue]{hyperref}
\usepackage{enumitem}
\usepackage[utf8]{inputenc}
\usepackage[english]{babel}
\usepackage[table,xcdraw]{xcolor}
\usepackage{booktabs}
\usepackage{authblk}
\usepackage{verbatim}
\usepackage{tikz}
\usepackage{pgfkeys}
\usetikzlibrary{fit}
\usetikzlibrary{shapes.geometric}
\usetikzlibrary{decorations.pathreplacing}
\usetikzlibrary{calc}

\usepackage{hyperref}

\newcommand{\Z}{\mathbb{Z}}

\newcommand{\DE}{\text{(DE)}}

\newcommand{\bsubeq}{\begin{subequations}}
\newcommand{\esubeq}{\end{subequations}}
\newcommand{\BI}{\begin{itemize}}
\newcommand{\EI}{\end{itemize}}

\newcommand{\cred}{\color{black}}
\newcommand{\ccred}{\color{black}}

\newcommand{\mc}{\mathcal}

\def\st{{\rm s.t.}}

\usepackage{algorithmic}
\usepackage{algorithm}
\usepackage{graphicx,amssymb}
\usepackage{enumerate}
\usepackage{subfigure}
\usepackage[normalem]{ulem}

\newcommand{\be}{\begin{enumerate}}
\newcommand{\ee}{\end{enumerate}}
\newtheorem{theorem}{Theorem}
\newtheorem{example}{Example}

\newcommand{\Azero}{\mc{A}^s_{hdr0}}
\newcommand{\Aone}{\mc{A}^s_{hdr1}}

\newcommand{\Aonep}{\mc{A}^s_{hdpr1}}

\usepackage{longtable}
\usepackage{afterpage}
\newcommand{\sm}{\setminus}

\usepackage{mathtools}

\newcommand{\Exp}{\mathbb{E}}

\usepackage{multirow}
\usepackage{graphicx}
\usepackage{caption}
\usepackage{breqn}
\allowdisplaybreaks

\usepackage{lipsum}
\makeatletter
\renewcommand\footnoterule{%
  \kern-3\p@
  \hrule\@width \textwidth
  \kern2.6\p@}
\makeatother

\makeatletter
\renewcommand*{\@fnsymbol}[1]{\ensuremath{\ifcase#1\or *\or \dagger\or \ddagger\or **\or
   \mathsection\or \mathparagraph\or \|\or  \dagger\dagger
   \or \ddagger\ddagger \else\@ctrerr\fi}}
\makeatother

\usepackage{natbib}
\setcitestyle{authoryear} 
\setcitestyle{notesep={}}

\usepackage{hyperref}

\title {Logic-based Benders Decomposition and Binary Decision Diagram Based Approaches for Stochastic Distributed Operating Room Scheduling}

\begin{document}

\author[1]{Cheng Guo\thanks{cguo@mie.utoronto.ca}}
\author[1]{Merve Bodur\thanks{bodur@mie.utoronto.ca}}
\author[1, 2, 3]{Dionne M. Aleman\thanks{aleman@mie.utoronto.ca}}
\author[4,5]{David R. Urbach\thanks{david.urbach@wchospital.ca}}

\affil[1]{\small Department of Mechanical and Industrial Engineering, University of Toronto, Toronto, Ontario M5S 3G8, Canada}
\affil[2]{\small Institute of Health Policy, Management and Evaluation, University of Toronto, Toronto, Ontario M5S 3E3, Canada}
\affil[3]{\small Techna Institute at University Health Network, Toronto, Ontario M5G 1P5, Canada}
\affil[4]{\small Department of Surgery, Women's College Hospital, Toronto, Ontario M5S 1B2, Canada}
\affil[5]{\small Department of Surgery, University of Toronto, Toronto, Ontario M5S 3E3, Canada}

\date{}
\maketitle

%

\begin{abstract}
The distributed operating room (OR) scheduling problem aims to find an assignment of surgeries to ORs across collaborating hospitals that share their waiting lists and ORs. We propose a stochastic extension of this problem where surgery durations are considered to be uncertain. In order to obtain solutions for the challenging stochastic model, we use sample average approximation, and develop two enhanced decomposition frameworks that use logic-based Benders (LBBD) optimality cuts and binary decision diagram based Benders cuts. Specifically, to the best of our knowledge, deriving LBBD optimality cuts in a stochastic programming context is new to the literature. Our computational experiments on a hospital dataset illustrate that the stochastic formulation generates robust schedules, and that our algorithms improve the computational efficiency. 
\end{abstract}



\section{Introduction}
According to a {\ccred recent} report by the Canadian Institute for Health Information, about 30\% of patients who need to receive hip replacement, knee replacement, and cataract surgeries wait longer 
than the recommended waiting time \citep{canadian2020wait}. As a result of such lengthy waiting times for medically necessary treatments, patients suffer from lost wages and reduced productivity, due to the effects of an untreated medical condition on mind and body \citep{esmail2019private}. The research on operating room (OR) scheduling studies how to allocate OR resources more strategically to improve its utilization. In this work, we study the \emph{stochastic distributed OR scheduling (SDORS) problem}, where several hospitals share their surgery waiting lists and ORs, and collectively schedule patients. We acknowledge the stochasticity in surgery procedures and derive schedules that remain robust under such an uncertainty. 

Studies from both hospital practice and mathematical modeling demonstrate that by utilizing distributed OR scheduling (DORS), ORs see an increased utilization rate \citep{magnussen2007centralized, wang2016discrete}. This is because when hospitals collaborate in waiting list management, they achieve a more balanced schedule across hospitals.

On the other hand, including multiple hospitals in the OR scheduling process adds to its uncertainty, such as the uncertainty in the surgery durations. This may result in an increased number of surgery cancellations, which can be costly for both the cancelled patients and the hospital. {\cred In this paper, we consider day of the surgery cancellations due to OR overtime. Lack of OR time can be a significant contributor to same-day cancellations; for instance \cite{dimitriadis2013challenge} found that lack of OR time caused 17.31\% of all same-day cancellations.}

In order to derive more robust schedules in the face of such a stochasticity, we use stochastic optimization. More specifically, we model SDORS as a two-stage stochastic integer program (2SIP), which is widely used in planning and scheduling problems. In our problem we have continuous random variables in the 2SIP, and such models are usually reformulated via sample average approximation (SAA). However, mostly due to integer decision variables in both stages, solving the SAA reformulation directly with commercial solvers can take a prohibitively long time, especially when the instance size grows. Furthermore, the presence of integer variables, particularly in the second stage, limits the choice of decomposition algorithms that can be employed, and highly impact the algorithmic efficiency. 

In this work, we develop two decomposition algorithms that are applicable when there are integer variables in the decomposed subproblems. Those decomposition algorithms use logic-based Benders (LBBD) cuts and Binary Decision Diagram (BDD) based Benders cuts. We also incorporate classical Benders cuts from the linear programming (LP) relaxations of subproblems, which helps to improve the convergence. Moreover, we propose several algorithmic enhancements, including adapted first fit decreasing (FFD) heuristics to find initial solutions, relaxations of the subproblems that are used to tighten the master problem, and an early stopping scheme to eliminate suboptimal master solutions faster, which significantly improve computational efficiency.

The contributions of this work lies both in the modeling and the solution methodology. In the modeling side, we propose the SDORS model, and evaluate its value against the DORS model {\cred from the literature \citep{roshanaei2017propagating} }. In the methodology aspect, we apply LBBD on a 2SIP and derive LBBD optimality cuts. To our best knowledge, the application of LBBD optimality cuts in stochastic programming is new to the literature. We also show the efficacy of the BDD-based decomposition algorithm, proposed by \cite{lozano2018binary}, for a new problem class, while in the literature it was tested only on the traveling salesman problem with time windows. Note that SDORS is a special type of planning and scheduling problem. Therefore, the decomposition algorithms we derive for SDORS are also applicable for other planning and scheduling models that share similar structures.

The rest of this paper is organized as follows: Section \ref{chapter: lit_review} provides a review of the literature on OR scheduling problems and 2SIP algorithms, {\cred also detailing some important concepts used in our work}. Section \ref{chapter: problem} presents a detailed description of the SDORS problem and its formulation. Section \ref{chapter: saa} explains the SAA framework. Section \ref{chapter: algorithm1} and Section \ref{chapter: algorithm2} introduce a two-stage and a three-stage decomposition algorithm, respectively, to solve the SAA problem. Finally, Section \ref{chapter: result} presents the experimental results.
\section{Literature Review}\label{chapter: lit_review}
The OR scheduling problem, which involves different types of planning decisions, studies the allocation of resources related to surgeries. {\cred In this work, we include both hospital and OR opening decisions as well as surgery scheduling decisions, under the uncertainty of surgery durations. Extensive literature reviews on OR scheduling can be found in \cite{cardoen2010operating} and \cite{guerriero2011operational}. We focus our review on objective and decision variable setups, and on stochastic OR scheduling models.}

The OR scheduling is well-studied in the {\it deterministic setting}, where we assume all parameters are known. {\cred A variety of objectives are used in the scheduling models. \cite{marques2012integer} propose a mixed-integer programming (MIP) model that both improves the OR occupancy rate and shortens the surgery waiting time. \cite{fei2008solving} reduce the underutilization and overtime of ORs. \cite{roshanaei2017propagating} minimize the hospital and OR opening costs, and try to schedule patients based on their priority scores. Other goals of OR scheduling include improving the throughput \citep{harper2002framework}, avoiding peaks in resource occupancy \citep{belien2008branch}, and incorporating preferences of different parties \citep{blake2002goal}. OR scheduling models in the literature also include different types of decisions. \cite{blake2002mount} assign ORs to different specialties in a hospital. \cite{belien2008branch} study the assignment of surgeons to ORs. \cite{jebali2006operating} consider both the assignment of patients to ORs and the sequencing of patiens in an OR. Similar to \cite{roshanaei2017propagating}, our model makes opening decisions and patient-to-OR assignment decisions, and our objective is to minimize the opening costs while accounting for patients' priority scores. In addition, we include cancellation decisions for overtime surgeries, and try to reduce surgery cancellations.} 


In real life, the scheduling of ORs can be affected by {\it uncertainty} in surgery procedures. By considering the uncertainty in OR scheduling, we can generate a more robust surgery schedule and achieve cost savings \citep{min2010scheduling}. {\cred Uncertainties in different stages of the process are considered in the literature. }\cite{denton2010optimal} minimize the cost in an OR scheduling problem assuming stochastic overtime. 
\cite{deng2019chance} deal with uncertainty in both waiting times and surgery durations. {\cred \cite{bowers2004managing} consider the uncertainty in patient arrivals. Also, stochastic OR scheduling problems in the literature use a variety of modeling techniques, including 2SIP \citep{denton2010optimal, min2010scheduling}, robust optimization \citep{denton2010optimal}, distributionally robust chance constraints \citep{deng2019chance}, and {\ccred multi-stage stochastic programming \citep{gul2015progressive}}. In our work, we formulate a 2SIP model which incorporates the stochasticity from surgery durations.}  

The importance of {\it distributed} planning and scheduling gained attention in many industries including manufacturing and healthcare. In a distributed setting, several agents share their resources and the list of tasks. Distributed manufacturing has the advantage of improving product quality, reducing costs and management risks, which helps the manufacturer to become more competitive under globalization \citep{naderi2014scatter}. \cite{timpe2000optimal} 
provide MIP formulations to solve lot-sizing problems in a distributed setting. 
\cite{behnamian2013heterogeneous} study a multi-factory production problem to minimize the maximum makespan. Inspired by the practice of hospitals in Toronto, \citet{roshanaei2017collaborative} propose a MIP model for the DORS problem where hospitals share their waiting lists and ORs. They find that in the distributed setting, hospitals are able to significantly reduce costs thanks to increased patient admission rate and OR throughput. For a similar setting, \citet{roshanaei2017propagating} propose an LBBD algorithm for DORS to solve practical instances efficiently. Our work is closely related to \citet{roshanaei2017propagating} and \citet{roshanaei2017collaborative}, as we also study the OR scheduling problem in a distributed setting. However, those previous works do not consider uncertainty in surgery durations, thus their proposed schedules may result in cancellations due to surgeries exceeding OR operating time limits. Moreover, our proposed SDORS model is significantly more challenging to solve, which calls for a thorough investigation of the problem structure to design efficient algorithms. 


SAA is the most common approach for 2SIPs, where the SAA problem, also known as the {\it deterministic equivalent} (DE), is usually solved {\cred exactly} by means of {\it decomposition} dividing the problem into a master problem and a subproblem. 
Benders decomposition is one of the most commonly-used decomposition algorithms \citep{benders1962partitioning}. However, it is not applicable when the subproblem contains  integer variables. \cite{laporte1993integer} propose the integer L-Shaped method for 2SIPs, which allows the subproblem to contain integer variables. \cite{angulo2016improving} present an improved integer L-Shaped method. {\cred 
\cite{caroe1998shaped} propose a generalized Benders decomposition method for 2SIP, which might generate nonlinear cuts and is computationally expensive in general cases. \cite{sherali2002modification} use the Reformulation-Linearization Technique to generate cuts from integer recourse problems. \cite{gade2014decomposition} use parametric Gomory cuts in their decomposition algorithm, which shows effectiveness for a general class of 2SIPs. Some other cutting planes for 2SIPs with discrete recourse include the lift-and-project cuts \citep{balas1993lift, ntaimo2008computations}, Fenchel decomposition (FD) cuts \citep{ntaimo2013fenchel}, and scenario FD cuts \citep{beier2015stage}. For a more detailed review on exact decomposition methods for 2SIPs, we refer the readers to \cite{bodur2015valid}. In this work, we instead choose to use LBBD and BDD-based methods, to benefit from our specific problem structure.}

The  \emph{LBBD method} proposed by \cite{hooker2003logic} provides a general way to decompose the problem such that the master problem and subproblem can be of any structure. {\cred Similar to the classical Benders method, the LBBD method generates LBBD feasibility cuts when the subproblem is infeasible, while LBBD optimality cuts are generated when the master problem underestimates the subproblem objective. However, both types of LBBD cuts are problem-specific, and their strength heavily depends on the incorporation of the underlying structure information. When strong cuts are derived, LBBD can be much more efficient than generic 2SIP methods.} In some recent works, LBBD is used to solve large-scale 2SIPs. 
\cite{lombardi2010stochastic} solve a stochastic allocation and scheduling problem for a multi-processor system-on-chips problem with LBBD. \cite{fazel2013solving} propose an LBBD method for a facility location/fleet management problem. {\cred Those papers} formulate {\it only LBBD feasibility cuts} {\cred in} the decomposition. In our work, we propose \emph{LBBD optimality cuts} for our 2SIP SDORS model in two different  decomposition frameworks. 


Our proposed decomposition algorithms for SDORS also use \emph{BDD-based cuts}. BDD is a graph structure that transforms a binary integer program (BIP) to an LP using a recursive formulation of the original problem. \cite{lozano2018binary} propose a decomposition algorithm based on BDD to solve a special class of 2SIPs, where 
each of the constraints linking the first and second stages consists of some binary second stage variables and a single first stage variable, and is deactivated when the corresponding first stage variable is zero. 
They transform the recourse problem to a capacitated shortest path problem, which has an LP formulation thus {\cred leads} to a  classical Benders decomposition algorithm. They conduct numerical experiments on the traveling salesman problem with time windows and show that their algorithm achieves substantial speedups compared to a commercial IP solver. In our work, we adapt an enhanced version of their decomposition algorithm to solve a 2SIP in the distributed OR scheduling setting.


\section{Problem Definition and Formulation}\label{chapter: problem}
In the SDORS problem, we aim to assign surgeries to the ORs in collaborative hospitals in the current planning horizon. The surgery durations are assumed to be stochastic. For each OR, there is a daily operating time limit. If the total surgery time of the day is expected to exceed the limit, some of the surgeries will be cancelled to satisfy this constraint. Our goal is to minimize total costs while ensuring more emergent patients get scheduled first.

We formulate the SDORS problem as a 2SIP. We assume that the probability distributions of surgery durations are known. In a 2SIP the decision process is divided into two stages. The first stage problem is solved before the revelation of surgery durations, where we make the decisions of (1) which {\cred surgical suite} to open during the current planning horizon, (2) which ORs to open in an opened hospital, (3) which patients to assign to the opened ORs, and (4) which patients to postpone to a future planning horizon. {\cred Notice that we consider only one list of patients in our model, for instance those patients may be from a single specialty, and we consider a surgical suite open if it is going to be (partially) used for our targeted patient list.}

{\cred In this paper, we consider day of the surgery cancellations. A 2013 study \citep{dimitriadis2013challenge} shows that the same-day surgery cancellation rate in the study is 5.19\%. Among all same-day cancellations, 17.31\% are due to lack of theatre time. We assume} each OR has an operating time limit 
for the total duration of the surgeries during a day, any surgery that is expected to finish after that time limit will be cancelled before it starts. The cancelled surgeries need to be scheduled in a future planning horizon instead. Therefore, each cancellation is associated with some costs. The second stage of the 2SIP is a recourse problem that takes the surgery-OR assignment decisions from the first stage, and decides on which surgeries to cancel, while minimizing the total cancellation cost.  {\cred Notice that 2SIP models implicitly assume that in the second stage, all surgery durations are known with certainty at the same time. In addition, in real life cancellation decisions rely on the sequence of surgeries, which we do not consider in our model. We note that the implementable decisions from our model are only the first-stage ones, i.e., opening and assignment decisions. By considering surgery durations and cancellations in different scenarios, we are able to create more flexible schedules than the deterministic models, and reduce the occurrences of cancellations. Incorporation of surgery sequencing will make our model more realistic and possibly yield better schedules. However, the joint optimization of our current decisions and the sequencing ones would become computationally very difficult. Therefore, we leave it for future research and propose our current model as a first step in showing the importance of incorporating surgery duration uncertainty to the {\it distributed} OR scheduling problems. We also observe that many stochastic scheduling models (in a variety of application domains) make the simplifying assumption to ignore sequencing decisions, such as \cite{fei2009solving}, \cite{gul2015progressive} and \cite{bodur2016mixed}.  In particular, similar to our paper, \citet{gul2015progressive} assume no sequencing and all surgery durations are revealed at the same time. Their results show decreased costs compared with deterministic models. Therefore, we believe our model is still helpful in improving the quality of distributed OR scheduling under uncertainty, despite having some simplifying assumptions. }

Except {\cred for} the stochasticity of surgery durations and cancellation decisions, all other settings in the SDORS model are the same as in {\cred the DORS model of \cite{roshanaei2017propagating}}. More specifically, for hospitals in the set $\mc{H}$ we need to decide which patients to schedule in the current planning horizon $\mc{D}$, which is a set of days. Each hospital $h$ has a set of ORs, $\mc{R}_h$. For simplicity we assume that all ORs in the set $\mc{R}_h$ are homogenous, and all hospitals have the same number of ORs, none of which is a limiting assumption. {\cred For example, in real life ORs could be nonidentical, which means some surgeries can only be practiced in a subset of all available ORs. This type of requirement can be modelled by extra constraints such as in \cite{roshanaei2017collaborative}.} Each patient in the set $\mc{P}$ has a health status score $\omega_p$.
Given a health score threshold $\Gamma$, all patients whose health score $\omega_p \geq \Gamma$ are called {\it mandatory patients} and denoted by the set $\mc{P}'$, and the rest of the patients are defined as {\it non-mandatory patients}. Mandatory patients have to be scheduled in the current planning horizon, while non-mandatory patients can be postponed to a future planning horizon. The objective is to minimize the total cost, which we will explain in more details later in this section.

The notations we use in the model are listed in Table \ref{table: SDORS_notation}.

{
\small
\begin{longtable}[htbp]{ll}
\caption{Notation}\label{table: SDORS_notation}\\
\hline
\multicolumn{2}{l}{\textbf{Sets:}}       \\
$\mc{P}$  & Set of patients, $p\in \mc{P}$ \\
$\mc{P}'$   & Set of mandatory patients \\ 
$\mc{H}$	& Set of hospitals, $h\in \mc{H}$\\
$\mc{D}$	& Set of days in the current planning horizon, $d\in \mc{D}$\\
$\mc{R}_h$ & Set of ORs in the surgical suite of hospital $h$, $r\in \mc{R}_h$\\
$\mc{S}$ & Set of possible scenarios of uncertain surgery durations, $s\in \mc{S}$\\
\multicolumn{2}{l}{\textbf{Parameters:}} \\
$G_{hd}$   & Cost of opening the surgical suite of hospital $h$ in day $d$\\
$F_{hd}$	&Cost of opening an OR in hospital $h$ on day $d$\\
$B_{hd}$	&Operating time limit of each OR on day d in hospital $h$\\
$T_p$	&Total booked time of patient $p$ \\ 
$c^{\text{sched}}_{dp}$ &Benefit for scheduled patient $p$ whose surgery is scheduled on day $d$\\
$c^{\text{unsched}}_{p}$ &Penalty for unscheduled patient $p$\\
$c^{\text{cancel}}_{p}$ &Penalty for cancelled patient $p$\\*[0.25cm]
 \multicolumn{2}{l}{\textbf{Decision variables:}} \\
$u_{hd}$ &1 if the surgical suite in hospital $h$ is opened on day $d$, 0 otherwise\\
$y_{hdr}$ &1 if OR $r$ of hospital $h$ is opened on day $d$, 0 otherwise\\
$x_{hdpr}$ &1 if patient $p$ is assigned to OR $r$ of hospital $h$ on day $d$, 0 otherwise\\
$w_p$  &1 if patient $p$ is postponed to a future planning horizon, 0 otherwise\\ 
$z_{hdpr}$ &1 if patient $p$'s surgery in OR $r$ of hospital $h$ on day $d$ is operated, 0 if it is cancelled \\ \hline
\end{longtable}
}

The SDORS problem is formulated as follows:
\begin{subequations}\label{model: sdors}
\begin{alignat}{5}
\text{(SDORS):} \ \min~~&\sum_{h\in\mc{H}} \sum_{d\in\mc{D}} G_{hd}u_{hd}+\sum_{h\in\mc{H}}\sum_{d\in\mc{D}}\sum_{r\in\mc{R}_h}F_{hd}y_{hdr}+\sum_{h\in\mc{H}}\sum_{d\in\mc{D}}\sum_{p\in\mc{P}}\sum_{r\in\mc{R}_h}c^{\text{sched}}_{dp}x_{hdpr}+  \hspace*{-5cm}\nonumber\\
&\sum_{p\in\mc{P} \setminus \mc{P'}}c^{\text{unsched}}_{p}w_p+\Exp_{\bf T} \mc{Q}({\bf x, y, T}) \hspace*{-5cm} &&\label{eq: sdors_obj}\\
\st ~~& \sum_{h\in\mc{H}}\sum_{d\in\mc{D}}\sum_{r\in\mc{R}_h}x_{hdpr}=1 & & \forall p \in \mc{P'} \label{eq: sors_assign_madatory_patients}\\
      & \sum_{h\in\mc{H}}\sum_{d\in\mc{D}}\sum_{r\in\mc{R}_h}x_{hdpr}+w_p=1 & & \forall p \in \mc{P}\setminus\mc{P'} \label{eq: sors_assign_nonmandatory_patients}\\
      & y_{hdr}\leq y_{hd,r-1} & &\forall h\in \mc{H}, d\in\mc{D},r\in\mc{R}_h\setminus\{1\} \label{eq: sors_symmetry_breaking_y} \\
       & \sum_{p\in\mc{P}}c^{\text{cancel}}_{p}x_{hdpr}\leq \sum_{p\in\mc{P}}c^{\text{cancel}}_{p}x_{hdp,r-1} & &\forall h\in \mc{H}, d\in\mc{D},r\in\mc{R}_h\setminus\{1\} \label{eq: sors_symmetry_breaking_x} \\
     & y_{hdr}\leq u_{hd}& &\forall h\in \mc{H}, d\in\mc{D}, r\in\mc{R}_h \label{eq: sors_y-u} \\
     & x_{hdpr}\leq y_{hdr} & &\forall h\in \mc{H}, d\in\mc{D}, p\in \mc{P}, r\in\mc{R}_h \label{eq: sors_x-y} \\
     & u_{hd},y_{hdr},x_{hdpr}, w_p\in \{0,1\} & &\forall h\in \mc{H}, d\in\mc{D}, p\in \mc{P}, r\in\mc{R}_h\label{eq: sors_binary}
\end{alignat}
\end{subequations}
where 
\begin{subequations}\label{model: sdors-sec}
\begin{alignat}{5}
\mc{Q}({\bf x, y, T}) = \min~~& \sum_{h\in\mc{H}}\sum_{d\in\mc{D}}\sum_{p\in\mc{P}}\sum_{r\in\mc{R}_h}c^{\text{cancel}}_p ({\cred x_{hdpr}} - z_{hdpr})  \label{eq: sdors_obj2}\\
\st ~~&\sum_{p\in\mc{P}} T_p z_{hdpr} \leq B_{hd}y_{hdr} &~~~& \forall h\in\mc{H}, d\in\mc{D}, r\in\mc{R}_h  \label{eq: sors_time_limit} \\
& z_{hdpr}\leq x_{hdpr}  & & \forall h\in\mc{H}, d\in\mc{D}, p\in\mc{P}, r\in\mc{R}_h \label{eq: sors_z-x} \\
&z_{hdpr}\in \{0,1\} & &\forall h\in \mc{H}, d\in\mc{D}, p\in \mc{P}, r\in\mc{R}_h
\end{alignat}
\end{subequations}
representing the optimization problem that minimizes the second stage cancellation cost, parameterized by the first stage decisions ${\bf x} = \{x_{hdpr} :  \forall h\in \mc{H}, d\in\mc{D}, p\in \mc{P}, r\in\mc{R}_h\}$, ${\bf y} = \{y_{hdr} : \forall h\in \mc{H}, d\in\mc{D}, r\in \mc{R}_h\}$, and ${\bf T} = \{T_1, ..., T_{|\mc{P}|}\}$. Note that ${\bf T}$ is the set of given parameters for surgery durations. 

The first stage objective function \eqref{eq: sdors_obj} minimizes the total operational and expected cancellation {\cred costs}, whose terms are respectively (i) the total cost of opening surgical suites in hospitals, (ii) the total cost of opening ORs, (iii) the benefit of scheduling patients (note that this term is negative), (iv) the penalty for postponing patients, and (v) the expected total cancellation cost. Constraints \eqref{eq: sors_assign_madatory_patients} ensure that mandatory patients are all scheduled in the current planning horizon. Constraints \eqref{eq: sors_assign_nonmandatory_patients} either schedule non-mandatory patients in the current planning horizon or postpone their surgeries to a future time. Constraints \eqref{eq: sors_symmetry_breaking_y} are symmetry breaking constraints for ORs in the same hospital, since they are homogeneous. Constraints \eqref{eq: sors_symmetry_breaking_x} are symmetry breaking constraints ensuring that only one permutation of patient assignment is allowed in each hospital-day pair $(h, d)$. {\cred More specifically, these constraints require that for ORs in the same hospital, the total cancellation cost of patients in an OR should be higher if this OR has a lower index. These constraints break more symmetry if cancellation costs are different for different patients, which is generally the case in our data. }Constraints \eqref{eq: sors_y-u} link $u$ and $y$ variables, to ensure that if the surgical suite in a hospital is closed, no OR in it is opened. Similarly, constraints \eqref{eq: sors_x-y} make sure that a patient will only be assigned to an opened OR. Lastly, constraints \eqref{eq: sors_binary} enforce binary restrictions on the decision variables.

For the second stage, i.e., the recourse problem \eqref{model: sdors-sec}, constraints \eqref{eq: sors_time_limit} ensure that the surgeries which are eventually conducted in an OR will not exceed its operating time limit. Constraints \eqref{eq: sors_z-x} link variables $z$ and $x$, which makes sure a patient's surgery can only be performed if it is scheduled in the first stage. Note that the second stage problem is always feasible, regardless of the first stage solution, because in the worst case we can cancel all the scheduled patients. Therefore, (SDORS) has {\it complete recourse}. 
%


\section{The SAA Problem}\label{chapter: saa}
Solving (SDORS) exactly involves calculating the multidimensional expectation over a set of continuous random vector ${\bf T}$. To overcome this difficulty, we solve the problem approximately by SAA. We generate a set of {\it scenarios} for the surgery duration vector ${\bf T}$. The set of scenarios are denoted by $\mc{S}$ and the surgery duration of patient $p$ under scenario $s\in\mc{S}$ is denoted by $T^s_p$. Under each scenario we need to decide whether to operate or cancel a surgery, thus we now have a decision variable $z^s_{hdpr}$ for each scenario that represents such decisions. We assume that all scenarios are equally likely, so each has a realization probability of $1 / |\mc{S}|$.

Replacing the expected cancellation cost with the average cancellation cost over the scenarios, we obtain the DE of the 2SIP:
\begin{subequations}\label{model: saa}
\begin{alignat}{5}
\DE : \min~~&\sum_{h\in\mc{H}} \sum_{d\in\mc{D}} G_{hd}u_{hd}+\sum_{h\in\mc{H}}\sum_{d\in\mc{D}}\sum_{r\in\mc{R}_h}F_{hd}y_{hdr}+\sum_{h\in\mc{H}}\sum_{d\in\mc{D}}\sum_{p\in\mc{P}}\sum_{r\in\mc{R}_h}c^{\text{sched}}_{dp}x_{hdpr}+  \hspace*{-6cm}\nonumber\\
&\sum_{p\in\mc{P} \setminus \mc{P'}}c^{\text{unsched}}_{p}w_p+ \frac{1}{|\mc{S}|} \sum_{s\in\mc{S}}\sum_{h\in\mc{H}}\sum_{d\in\mc{D}}\sum_{p\in\mc{P}}\sum_{r\in\mc{R}_h}c^{\text{cancel}}_p (x_{hdpr} - z^s_{hdpr})  \label{eq: saa_obj}\hspace*{-5cm} &&\\
\st ~~&\eqref{eq: sors_assign_madatory_patients} - \eqref{eq: sors_binary}\\
       & \sum_{p\in\mc{P}} T^s_p z^s_{hdpr} \leq B_{hd}y_{hdr} && \forall h\in\mc{H}, d\in\mc{D}, r\in\mc{R}_h, s\in\mc{S} \label{eq: saa_time_limit} \\
       & z^s_{hdpr}\leq x_{hdpr}  & & \forall h\in\mc{H}, d\in\mc{D}, p\in\mc{P}, r\in\mc{R}_h,  s\in\mc{S} \label{eq: saa_z-x} \\
     & z^s_{hdpr}\in \{0,1\} & &\forall h\in \mc{H}, d\in\mc{D}, p\in \mc{P}, r\in\mc{R}_h, s\in\mc{S}
\end{alignat}
\end{subequations}
Next, we present two decomposition schemes to solve this model exactly.
\section{Two-stage Decomposition}\label{chapter: algorithm1}
In this section we first introduce the decomposition scheme that divides the \DE~into a master problem and a set of subproblems, where each subproblem is a BIP. We generate three types of cutting planes from the subproblems: LBBD optimality cuts, classical Benders cuts from their BDD reformulations (BDD-based Benders cuts), and classical Benders cuts from their LP relaxations. Notice that the first two cut families are exact, meaning solely adding violated cuts from one such family will lead the algorithm to converge to an optimal \DE~solution. We also introduce several algorithmic enhancements for the decomposition, some of which are novel to the literature. 
\subsection{Decomposition Framework}
We decompose our problem into a master problem and a set of subproblems. In the master problem the first-stage decisions, namely hospital opening, OR opening, patient assignment, and postponing decisions, are made. In the subproblems we make decisions about surgery cancellation. More specifically, the master problem contains the variables $u_{hd}$, $y_{hdr}$, $x_{hdpr}$, and $w_p$ from the original problem, and a new variable $Q^s_{hdr}$ representing the cancellation cost of the hospital-day-OR combination $(h, d, r)$ under scenario $s$: %
 \begin{subequations}\label{model: sors_master}
\begin{alignat}{5}
\min~~&\sum_{h\in\mc{H}} \sum_{d\in\mc{D}} G_{hd}u_{hd}+\sum_{h\in\mc{H}}\sum_{d\in\mc{D}}\sum_{r\in\mc{R}}F_{hdr}y_{hdr}+\sum_{h\in\mc{H}}\sum_{d\in\mc{D}}\sum_{p\in\mc{P}}\sum_{r\in\mc{R}}{\cred c^{\text{sched}}_{dp}} x_{hdpr}\hspace{-3cm}\nonumber\\
&+ \sum_{p\in\mc{P} \setminus \mc{P'}}c^{\text{unsched}}_p w_p+ \frac{1}{|\mc{S}|}\sum_{s\in\mc{S}}\sum_{h\in\mc{H}}\sum_{d\in\mc{D}}\sum_{r\in\mc{R}} Q^s_{hdr} \label{eq: sors_master_obj}\\
\st ~~&\eqref{eq: sors_assign_madatory_patients} - \eqref{eq: sors_binary}\\
         & [\text{LBBD cuts or BDD-based Benders cuts}] \label{eq: sors_master_cuts}\\
         &Q^s_{hdr}\geq 0, &&\forall s\in\mc{S}, h\in\mc{H}, d\in\mc{D}, r\in\mc{R}_h
\end{alignat}
\end{subequations}

Constraints \eqref{eq: sors_master_cuts} are BDD-based Benders cuts or LBBD cuts that are generated progressively from the subproblems, in order to correctly approximate the second-stage value function, which is the link between the $Q$ variables and the first-stage variables. Without cutting planes \eqref{eq: sors_master_cuts}, the master problem \eqref{model: sors_master} is a relaxation of \DE. The value of cancellation cost $Q^s_{hdr}$ is underestimated because the master problem lacks enough information regarding the recourse decisions. On the other hand, subproblems contain the information about the recourse decisions. We provide such information by generating cutting planes from subproblems and adding them to the master problem. 

By solving the master problem, we get the optimal solutions for $u_{hd}$, $y_{hdr}$, $x_{hdpr}$, $w_p$ and $Q^s_{hdr}$, denoting their optimal solutions respectively by $\hat{u}_{hd}$, $\hat{y}_{hdr}$, $\hat{x}_{hdpr}$, $\hat{w}_p$, and $\hat{Q}^s_{hdr}$. Then we pass optimal solutions of $x_{hdpr}$ and $y_{hdr}$ to the subproblems. Note that we formulate a subproblem for each $(h, d, r, s)$ tuple, as the recourse decisions for each OR and scenario is independent once the values of $x_{hdpr}$ and $y_{hdr}$ are fixed. Each subproblem minimizes the cancellation cost $Q^s_{hdr}$ by selecting the least costly patients to cancel, if the current assignment exceeds the operating time limit of an OR:
\begin{subequations}\label{model: sors_sub}
\begin{alignat}{5}
{\cred \mc{Q}^s_{hdr}(\hat{x}_{hd\cdot r}, \hat{y}_{hdr}, T^s_\cdot)}= \min~~& \sum_{p\in \mc{P}} c^{\text{cancel}}_p (\hat{x}_{hdpr} - z^s_{hdpr})\hspace*{3cm}&&& \label{eq: sors_sub_obj}\\
\st ~~& \sum_{p\in \mc{P}}T^s_p z^s_{hdpr} \leq B_{hd} && \label{eq: sors_sub_time_limit}\\
	& z^s_{hdpr} \leq \hat{x}_{hdpr} & & \forall p\in \mc{P} \label{eq: sors_sub_z-x}\\
     	& z^s_{hdpr}\in \{0,1\}& &\forall p\in \mc{P} \label{eq: sors_sub_binary}
\end{alignat}
\end{subequations}
{\ccred Note that we use $\cdot$ to represent all members of the index in that position, thus $\hat{x}_{hd\cdot r} = \{\hat{x}_{hd1r}, ..., \hat{x}_{hd|\mc{P}| r}\}$ and $T^s_\cdot = \{T^s_1, ..., T^s_{|\mc{P}|}\}$.} Constraints \eqref{eq: sors_sub_time_limit} are the time limit constraints and constraints \eqref{eq: sors_sub_z-x} guarantee that no surgery will be performed if it is not scheduled in the first place. {\cred Constraints \eqref{eq: sors_sub_binary} enforce variables $z^s_{hdpr}$ to be binary. To simplify the notation, from now on we use $\bar{Q}^s_{hdr}$ to denote $\mc{Q}^s_{hdr}(\hat{x}_{hd\cdot r}, \hat{y}_{hdr}, T^s_\cdot)$, which is the correct optimal value of the subproblem objective.}
\subsection{The LBBD cuts}
The first category of cutting planes that we derive are the LBBD cuts. Like the case with classical Benders cuts, there are two types of LBBD cuts: LBBD feasibility cuts and LBBD optimality cuts. LBBD feasibility cuts are generated when the master problem solution makes a subproblem infeasible. On the other hand, LBBD optimality cuts are generated when the master problem solution {\cred underestimates the (minimization) objective value of a feasible subproblem, }due to lacking the information about the subproblem. In our decomposition framework, the subproblems are never infeasible, since it can always cancel all the assigned patients and get a feasible solution. As such, we do not need LBBD feasibility cuts. 

Since the master problem is a relaxation of \DE, the cancellation cost $\hat{Q}^s_{hdr}$ is an underestimate of the optimal cancellation cost for the assignment $(\hat{u}_{hd}, \hat{y}_{hdr}, \hat{x}_{hdpr}, \hat{w}_p)$ from the master problem, while the subproblem objective $\bar{Q}^s_{hdr}$ provides the actual cancellation cost of such an assignment. When $\hat{Q}^s_{hdr} < \bar{Q}^s_{hdr}$, we add the following LBBD optimality cut to the master problem:
\begin{align}\label{eq: lbbd_cut}
{\cred Q^s_{hdr}} \geq \bar{Q}^s_{hdr} -\sum_{p \in \hat{\mc{P}}_{hdr}} c_p^{\text{cancel}}(1 - x_{hdpr})
\end{align}
where $\hat{\mc{P}}_{hdr}$ is the patient list corresponding to the master solution: $\hat{\mc{P}}_{hdr} = \{p\in\mc{P}~|~\hat{x}_{hdpr} = 1\}$. This cut provides a lower bound (LB) for $Q^s_{hdr}$: the cancellation cost corresponding to patients in $\hat{\mc{P}}_{hdr}$ is $\bar{Q}^s_{hdr}$, and if we remove a patient $p\in\hat{\mc{P}}_{hdr}$ from the current assignment, we can potentially save the corresponding cancellation cost, $c^{\text{cancel}}_p$.

Note that the LBBD cut is tight at the current master solution, i.e., it provides the exact cancellation cost for the master solution $(\hat{u}_{hd}, \hat{y}_{hdr}, \hat{x}_{hdpr}, \hat{w}_p)$. More importantly, different than the vanilla no-good cut, it also provides a (not necessarily tight) LB for any assignment with a patient list that is either a superset of $\hat{\mc{P}}_{hdr}$, or contains a subset of $\hat{\mc{P}}_{hdr}$. 

Theorem \ref{theorem: lbbd} states that the proposed LBBD cut is valid, which means it eliminates the current master solution, and more importantly it does not exclude any optimal solution of (DE). Considering that there are finitely many binary $x_{hdpr}$ feasible solutions in the master problem, this result implies that the LBBD algorithm with the cuts \eqref{eq: lbbd_cut} converge to an optimal solution of (DE) in a finite number of iterations.
\begin{theorem}\label{theorem: lbbd}
The LBBD optimality cut \eqref{eq: lbbd_cut} is valid.
\end{theorem}
The proof for the validity of the LBBD cuts is given in the Online Supplement. 

\subsection{BDD-based Benders Cuts}\label{chapter: bdd}
The next category of cutting planes, which we call {\it BDD-based Benders cuts}, are based on the work of \cite{lozano2018binary}. {\cred Those cuts are applicable here because (i) in our decomposition the master variables that are linked to the subproblems, {\ccred i.e. the x-variables}, are all binary; and (ii) any subproblem constraint that is impacted by the first-stage decisions is deactivated by a {\it single} $x_{hdpr}$ variable. Despite the fact that the subproblems are BIPs, the BDD-based Benders cuts are actually a set of classical Benders cuts; they are obtained by first transforming a subproblem into a shortest path problem, then generating classical Benders optimality cuts from the reformulation which is now an LP. In what follows, we first provide some basic concepts for BDDs, and then explain how to transform our subproblems into shortest path problems via BDDs and in turn obtain the additional set of Benders cuts.} 

{\cred A BDD is a layered directed acyclic graph with a single root node and a single terminal node, which is used to represent a BIP, e.g., see the example in Figure \ref{fig: bdd_example}. The arcs of a BDD correspond to assigning values to binary variables. More specifically, each layer of arcs corresponds to a binary variable. There are two types of arcs, namely zero-arc and one-arc, respectively corresponding to the associated variable taking the value of zero or one. There is a one-to-one correspondence between feasible solutions of the BIP and root-to-terminal paths in the BDD, thus the BDD compactly represents the feasible set of the BIP. Moreover, arcs in the BDD are assigned length values in such a way that for each root-to-terminal path, its path length is equal to the BIP objective function value of the corresponding BIP feasible solution. Therefore, finding an optimal solution to BIP with a minimization/maximization objective reduces to finding  a shortest/longest path in the BDD. For more details on BDDs, we refer the reader to the book by \cite{bergman2016decision}.}


In our problem, we transfer each {\cred subproblem \eqref{model: sors_sub} to a BDD. We note that although the BDD size can be exponential in the BIP size, it is pseudo-polynomial in our case since our subproblem is a knapsack problem. We first provide a simple example to illustrate the BDD transformation. 

\begin{example}
{\normalfont Suppose in {\cred the subproblem \eqref{model: sors_sub}} for an $(h, d, r, s)$ tuple there are four scheduled patients $(p = 1,...,4)$, and we need to decide whether to cancel them. {\ccred Assume that for those four patients their cancellation costs are respectively $4,1,3,8$ and their surgery durations in the current scenario are respectively $2,1,3,3$, then the} subproblem formulation is the following {\cred (we have replaced all ${\ccred \hat{x}_{hdpr}}, p =1,...,4$ in \eqref{model: sors_sub} with 1 as all four patients are selected)}:
\begin{subequations}\label{model: bdd_small_exam}
\begin{alignat}{4}
{\cred \bar{Q}_{hdr}^s = \min} \ & 4(1 - z^s_{hd1r})+(1- z^s_{hd2r})+3(1- z^s_{hd3r})+8(1-z^s_{hd4r})\label{eq: bdd_small_obj1} &&\\
\text{s.t.} \ & {\cred 2z^s_{hd1r}+z^s_{hd2r}+3z^s_{hd3r}+3z^s_{hd4r} \leq 5}&& \\
&z^s_{hdpr}  \leq 1 &&\forall p\in\hat{\mc{P}}_{hdr}\\
& z^s_{hdpr} \in \{0,1\} &&\forall p\in\hat{\mc{P}}_{hdr}
\end{alignat}
\end{subequations}
where $\hat{\mc{P}}_{hdr} = \{1,2,3,4\}$ {\cred represents the subset of selected patients. To simplify the notation, we equivalently rewrite the objective as:} 
\begin{align}\label{eq: bdd_small_obj2}
{\cred \bar{Q}_{hdr}^s - 16 = \min~~ - 4z^s_{hd1r} - z^s_{hd2r} - 3z^s_{hd3r} - 8z^s_{hd4r}}
\end{align}

The problem is represented via the BDD in Figure \ref{fig: bdd_example}. 
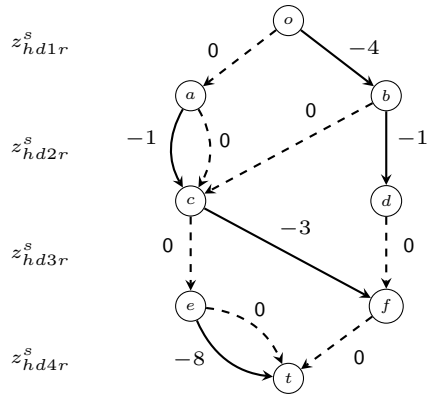
\begin{figure}[h]
\centering
\begin{tikzpicture}[scale=0.3][font=\sffamily,\tiny]
	\tikzset{mycircle/.style={draw,circle, minimum size=0.4cm}}
	\tikzstyle{arrow} = [thick,->,>=stealth]

	\node [mycircle] (r) at (0,0) {$o$};

	\node [mycircle,left of=r,xshift=-0.3cm,yshift=-1cm] (a)  {$a$};
	\node [mycircle,right of=r,xshift=0.3cm,yshift=-1cm] (b) {$b$};

	\node [mycircle,below of=a,yshift=-0.4cm] (c)  {$c$};
	\node [mycircle,below of=b,yshift=-0.4cm] (d) {$d$};

	\node [mycircle,below of=c,yshift=-0.4cm] (e)  {$e$};
	\node [mycircle,below of=d,yshift=-0.4cm] (f) {$f$};

	\node [mycircle] at ([yshift=-3.2cm]e -| r) (t) {$t$} ;

	\node [left of=a,xshift=-1cm,yshift=0.7cm,scale=1.2] (x1) {$z^s_{hd1r}$};
	\node [left of=a,xshift=-1cm,yshift=-0.7cm,scale=1.2] (x2) {$z^s_{hd2r}$};
	\node [left of=c,xshift=-1cm,yshift=-0.7cm,scale=1.2] (x3) {$z^s_{hd3r}$};
	\node [left of=e,xshift=-1cm,yshift=-0.7cm,scale=1.2] (x4) {$z^s_{hd4r}$};

	\draw[arrow,dashed] (r) -- node[anchor=south,xshift = -0.35cm,yshift=-0.1cm,scale=1.3] {0} (a);
	\draw[arrow] (r) -- node[anchor=south,xshift = 0.35cm,yshift=-0.1cm,scale=1.3] {$- 4$} (b);

	\draw[arrow,dashed] (a) [bend left] to node[anchor=south,xshift = 0.2cm,yshift=-0.1cm,scale=1.3] {0} (c);
	\draw[arrow] (a) [bend right] to node[anchor=south,xshift = -0.4cm,yshift=-0.1cm,scale=1.3] {$- 1$} (c);

	\draw[arrow,dashed] (b) -- node[anchor=south,xshift = 0.3cm,yshift=0.3cm,scale=1.3] {0} (c);
	\draw[arrow] (b) -- node[anchor=south,xshift = 0.35cm,yshift=-0.1cm,scale=1.3] {$- 1$} (d);
	\draw[arrow] (c) -- node[anchor=south,xshift = 0.1cm,yshift=0.1cm,scale=1.3] {$- 3$} (f);
	\draw[arrow,dashed] (c) -- node[anchor=south,xshift = -0.3cm,yshift=-0.1cm,scale=1.3] {0} (e);
	\draw[arrow,dashed] (d) -- node[anchor=south,xshift = 0.3cm,yshift=-0.1cm,scale=1.3] {0} (f);
	\draw[arrow,dashed] (f) -- node[anchor=south,xshift = 0.3cm,yshift=-0.4cm,scale=1.3] {0}(t);

	\draw[arrow,dashed] (e) [bend left] to node[anchor=south,xshift = 0.1cm,yshift=0cm,scale=1.3] {0} (t);
	\draw[arrow] (e) [bend right] to node[anchor=south,xshift = -0.5cm,yshift=-0.23cm,scale=1.3] {$- 8$} (t);
\end{tikzpicture}
\caption{BDD for a simple subproblem with four patients}
\label{fig: bdd_example}
\end{figure}
The BDD contains a root node at the top, denoted by $o$, and a terminal node at the bottom, denoted by $t$. Each layer of the arcs represents the values of the variables $z^s_{hdpr}$, from patient 1 to patient 4. (Note that different variable orderings can produce different BDDs, thus may impact the computational performance.) At each layer, the dashed and solid arcs respectively indicate selecting 0 and 1 values for the layer's associated variable. The value above or beside each arc is the cost of $z^s_{hdpr}$ if it receives the value indicated by the arc type. When transforming the subproblem to a shortest path problem, those values are used as costs for the arcs. Each path from the root node $o$ to the terminal node $t$ represents a feasible solution for the subproblem, and each feasible solution is represented by a path from  $o$ to $t$. Also, the path length provides the correct evaluation of the objective value for the corresponding feasible solution. {\cred For example, the path $o$-$b$-$d$-$f$-$t$ represents the feasible solution $z_{hd1r} = 1, z_{hd2r} = 1, z_{hd3r} = 0, z_{hd4r} = 0$. The length of this path is -5, which equals the objective value of the corresponding $z$ solution.} Therefore, by finding a shortest path from the root node to the terminal node, we can obtain an optimal solution to the subproblem.  \qed}

\end{example}
}
{\cred We now formally illustrate the BDD transformation of the subproblem \eqref{model: sors_sub}.} The following discussion deals with the subproblem corresponding to an $(h, d, r, s)$ tuple, and makes cancellation decisions concerning the {\cred scheduled} patient set $\hat{\mc{P}}_{hdr}$. We only consider patients from $\hat{\mc{P}}_{hdr}$ in the BDD reformulation. Let $\mc{G}^{s}_{hdr} = \{\mc{N}^s_{hdr},\mc{A}^s_{hdr}\}$ be the BDD for this subproblem where $\mc{N}^s_{hdr}$ is the set of nodes and $\mc{A}^{s}_{hdr}$ is the set of arcs. Let $\mc{A}^s_{hdr0}$ and $\mc{A}^s_{hdr1}$ be the sets of zero-arcs and one-arcs, i.e., corresponding to setting variables in $\{z^s_{hdpr}\}_{p\in\hat{\mc{P}}_{hdr}}$ equal to 0 and 1 in the subproblem, respectively. In addition, we introduce the sets $\mc{A}^s_{hdpr0}\subset\Azero$ and $\mc{A}^s_{hdpr1}\subset\Aone$ to respectively denote the sets of zero-arcs and one-arcs corresponding to a patient $p$. Nodes $o\in\mc{N}^s_{hdr}$ and $t\in\mc{N}^s_{hdr}$ are the root and terminal nodes of the diagram $\mc{G}^s_{hdr}$. The capacity of all arcs are one. We denote the cost of an arc $a\in \mc{A}^s_{hdr}$ as $g^{s}_{hdra}$ and its value is decided as follows: 
\[
    g^s_{hdra} = 
\begin{cases}
    - c^{\text{cancel}}_{p},& \text{if } a\in\Aonep \\
    0,              &\text{if } a\in\mc{A}^{s}_{hdr0}
\end{cases}
\]


Let the start and end points of arc $a$ be $s(a)$ and $d(a)$, respectively. Then the BDD reformulation of the subproblem as a capacitated shortest path problem is as follows:
\begin{subequations}\label{model: bdd_sub_trans}
	\begin{alignat}{5}
\hspace{-6mm}	\bar{Q}^s_{hdr} - \sum_{p\in\mc{P}}c^{\text{cancel}}_p \hat{x}_{hdpr} = \min~&\sum_{a \in \mc{A}^s_{hdr}} g^s_{hdra} f_a\\
	\st~&\sum_{a\in\mc{A}^s_{hdr} | s(a)=o} f_a=1&& &&&(\pi_o) \label{model: bdd_sub_trans_1}\\
	&\sum_{a\in\mc{A}^s_{hdr}|s(a)=i} f_a-\sum_{a\in\mc{A}^s_{hdr} | d(a)=i} f_a=0&~&\forall i\in\mc{N}^s\sm\{o,t\}&~~&&(\pi_i)\label{model: bdd_sub_trans_2}\\
	&\sum_{a\in\mc{A}^s_{hdr} | d(a) = t} f_a={\cred 1} && &&&(\pi_t)\label{model: bdd_sub_trans_3}\\
	&f_a\leq \hat{x}_{hdpr}  &&\forall a\in\mc{A}^{s}_{hdpr1}, p\in\mc{P}&&&(-\xi_a)\label{eq: bdd_sub_trans_fx}\\
	&f_a\geq 0 &&\forall a\in\mc{A}^s_{hdr}
	\end{alignat}
\end{subequations}
where $f_a$ is the variable for the flow through arc $a$. The variable in parentheses at the end of each constraint is the corresponding LP dual variable. Note that the objective value of problem \eqref{model: bdd_sub_trans} is different from that of the subproblem \eqref{model: sors_sub} by a constant term $\sum_{p\in\mc{P}}c^{\text{cancel}}_p \hat{x}_{hdpr}$. {\cred This is because in the process of the BDD transformation, we remove the constant term in \eqref{eq: sors_sub_obj}, $\sum_{p\in\mc{P}}c^{\text{cancel}}_p \hat{x}_{hdpr}$, just as shown in the small example, where the objective \eqref{eq: bdd_small_obj1} becomes \eqref{eq: bdd_small_obj2}.}


{\cred Let $\bar{\pi}_o$, $\bar{\pi}_i$ and $\bar{\pi}_t$ be the respective optimal dual values corresponding to constraints \eqref{model: bdd_sub_trans_1}, \eqref{model: bdd_sub_trans_2} and \eqref{model: bdd_sub_trans_3}. Also define the optimal value for the dual variable $\xi_a$ as $\bar{\xi}_a$.} 
We can derive the following classical Benders cut from the problem \eqref{model: bdd_sub_trans}, which we refer to as {\it BDD-based Benders cut}:
\begin{align}
Q^s_{hdr} - \sum_{p\in\mc{P}}c^{\text{cancel}}_p x_{hdpr} \geq {\cred \bar{\pi}_o} - \sum_{p\in\mc{P}}\Big(\sum_{a\in\mc{A}^{s}_{hdpr1}} \bar{\xi}_a\Big)x_{hdpr}
\end{align}
{\cred Note that $\bar{\pi}_t$ is set to 0 and thus not included in the BDD-based Benders cut. We are able to do this because the corresponding constraint of ${\pi}_t$, \eqref{model: bdd_sub_trans_3}, is linearly dependent on the other flow balance constraints \eqref{model: bdd_sub_trans_1}-\eqref{model: bdd_sub_trans_2}.}

As suggested by \cite{lozano2018binary}, this cut can be further strengthened by replacing the $x_{hdpr}$ coefficient of $\sum_{a\in\Aonep} \bar{\xi}_a$, with $\max_{a\in\Aonep}\bar{\xi}_a$: 
\begin{align}\label{eq: bdd_cut}
Q^s_{hdr} - \sum_{p\in\mc{P}}c^{\text{cancel}}_p x_{hdpr} \geq {\cred \bar{\pi}_o} -\sum_{p\in\mc{P}}\Big(\max_{a\in\mc{A}^{s}_{hdpr1}} \bar{\xi}_a\Big)x_{hdpr}
\end{align} 
In our implementation we use this strengthened version of BDD-based Benders cut. For the proof of the validity of \eqref{eq: bdd_cut}, we  refer the readers to \cite{lozano2018binary}. 
\subsection{Classical Benders Cuts from the Relaxed Subproblem}\label{chapter: benders}
In order to provide  a better LB for $Q^s_{hdr}$, we also add classical Benders cuts from the LP relaxation of the subproblem \eqref{model: sors_sub} to the master problem. Note that unlike the LBBD cuts and the BDD-based Benders cuts, the classical Benders cuts generated in this method cannot guarantee convergence. However, it speeds up the convergence in our algorithm when used together with LBBD cuts or BDD-based Benders cuts.

The process to obtain classical Benders cuts is standard. We linearly relax the subproblem \eqref{model: sors_sub} and denote the optimal objective value of the relaxed problem as $\bar{Q}_{hdr}^{s\text{LP}}$. Let 
$\bar{\delta}$ and $\bar{\eta}$ denote the optimal dual solutions that respectively correspond to constraints \eqref{eq: sors_sub_time_limit} and \eqref{eq: sors_sub_z-x} in the LP relaxation.
The classical Benders cut \eqref{eq: bd_cut}, which is generated and added to the master problem \eqref{model: sors_master} when $\hat{Q}^s_{hdr} < \bar{Q}^{s \text{LP}}_{hdr}$, is as follows:
\begin{align}\label{eq: bd_cut}
Q^s_{hdr} \geq \sum_{p\in \hat{\mc{P}}_{hdr}} (c_p^{\text{cancel}} + \bar{\delta}_p) x_{hdpr} + \bar{\eta} B_{hd}
\end{align}
\subsection{Additional Algorithmic Enhancements}\label{chapter: algo_enhancements}
In order to further improve the performance of the algorithm, we introduce three enhancement techniques. Section \ref{chapter: ffd} provides a heuristic way to find an initial solution, thus a good upper bound (UB) to start the search. Section \ref{chapter: qlb} derives a relaxation of the subproblem which is used to tighten the master formulation. Section \ref{chapter: heurcb} introduces a scheme that inserts additional heuristic solutions in the branch-and-cut algorithm to improve UBs. 
\subsubsection{Adapted First Fit Decreasing (FFD) Heuristic.}\label{chapter: ffd}
FFD is a heuristic method that finds a good feasible solution for the binary bin-packing problem \citep{johnson1974worst}. {\cred We choose to use this heuristic because it is quick to implement, and is flexible to adapt to our problem setup. More importantly, FFD used together with the LBBD method achieves good results in the literature of scheduling and assignment problems \citep{fazel2013solving, roshanaei2017propagating}.} We propose an adapted version of FFD to find an initial assignment of the patients for the SDORS problem. Viewing the ORs as bins and their operating time limits as the bin capacities, and the surgery duration of a patient as the weight of an item that needs to be fitted, the goal is to schedule as many patients to the ORs as possible. 

To find an initial assignment via the heuristic, we first sort the $(h,d)$ pairs in the decreasing order of the operating time limit $B_{hd}$, then sort the ORs from different $(h,d)$ pairs according to the sorted order of $(h,d)$ pairs. For ORs of the same $(h,d)$ pair, arrange them in the increasing order of their indices. Note that because we have both mandatory and non-mandatory patients, the FFD algorithm is adapted to accommodate the requirement that mandatory patients must be arranged in the current planning horizon. More specifically, when we initially sort the patients, we first sort all the mandatory patients, i.e. before sorting any non-mandatory patients. Both the mandatory and non-mandatory patients are sorted in the decreasing order of their health status score $\omega_p$.

Overall, the heuristic has the following steps: (i) Randomly pick a scenario $s\in\mc{S}$ (in the implementation we pick the first scenario). (ii) For each patient in the sorted order, check whether one of the opened ORs has enough capacity left to fit the patient based on the surgery duration in the selected scenario, starting from the first OR in the sorted order. If one of the ORs has enough capacity, assign the patient to the OR and update the capacity of the OR. If not, open a new OR to fit the patient in. (iii) If there are no new OR left to fit the patient, postpone the patient if she is not mandatory, otherwise {\cred the FFD heuristic cannot schedule all mandatory patients, and other heuristics could be tried. (In our experiments, we get feasible solutions from FFD for all instances, i.e., never encountered any instance where FFD cannot schedule mandatory patients.) }(iv) After assigning all the patients to ORs, we solve the subproblem \eqref{model: sors_sub} for each opened ORs to find its corresponding cancellation cost under this assignment. 

We use this heuristic solution as a warm start which provides a good initial UB for the algorithm. In addition, based on the heuristic solution, we derive LBBD cuts \eqref{eq: lbbd_cut} or BDD-based Benders cuts \eqref{eq: bdd_cut}, plus classical Benders cuts \eqref{eq: bd_cut}.   We add those cuts to the master problem \eqref{model: sors_master} before the branch-and-cut algorithm, which we will describe in more detail in Section \ref{subch: overall_implement}. Adding those cuts helps to reduce the initial optimality gap of the decomposition algorithm.
\subsubsection{Adding Subproblem Relaxations to the Master Problem.}\label{chapter: qlb}
In decomposition algorithms, it is observed that adding some form of relaxed subproblem to the master problem may greatly improve computational efficiency \citep{cire2016logic}. In our algorithm, we tighten the master problem by providing a LB for each $Q^s_{hdr}$, which can be derived by relaxing the subproblem:
\begin{align}\label{eq: qlb}
Q^s_{hdr} \geq \Bigg(\min_{p\in\mc{P}} \frac{c_p^{\text{cancel}}}{T^s_p}\Bigg) \Bigg(\sum_{p\in\mc{P}} T^s_p x_{hdpr} - B_{hd}\Bigg) ~~\forall h\in\mc{H}, d\in\mc{D}, r\in\mc{R}_h, s\in\mc{S}
\end{align}

Intuitively, those constraints approximate the OR cancellation cost by estimating the total overtime and the cost of cancelling overtime patients. We formally state the validity of those constraints in Theorem \ref{theorem: qlb} whose proof is provided in the Online Supplement. 

\begin{theorem}\label{theorem: qlb}
Constraints \eqref{eq: qlb} provide valid LBs for $Q^s_{hdr}$.
\end{theorem}

\subsubsection{Insert Heuristic Solution in Branch-and-Cut.}\label{chapter: heurcb}
{\cred While implementing the branch-and-cut algorithm, at each integral node we obtain integral solutions for the master variables. However, if the subproblem objectives are underestimated at those nodes, those master variable solutions are not utilized to obtain a UB. Inspired by the practice, e.g., in \cite{bodur2016mixed}, we combine the master variable solution at an integral node and the correct value for the corresponding subproblem, and provide this additional feasible solution to the solver. More specifically, each} additional solution is obtained at a master solution, $(\hat{u}, \hat{y}, \hat{x}, \hat{w}, \hat{Q})$, when $(\hat{u}, \hat{y}, \hat{x}, \hat{w})$ is {\it integral}. We solve the corresponding subproblems for this master solution and get $\bar{Q}$, then the solution $(\hat{u}, \hat{y}, \hat{x}, \hat{w}, \bar{Q})$ is a feasible solution for \DE. Next, we evaluate the objective value at $(\hat{u}, \hat{y}, \hat{x}, \hat{w}, \bar{Q})$. If it is better than the incumbent UB of branch-and-cut algorithm at the moment, then we update the incumbent UB by setting $(\hat{u}, \hat{y}, \hat{x}, \hat{w}, \bar{Q})$ as a heuristic solution via the commercial solver callback. %
{\ccred We} add additional solutions through the commercial solver callback inside the branch-and-cut algorithm. 
\subsection{Overall Implementation Approach}\label{subch: overall_implement}
As mentioned before, we add the cutting planes, including the LBBD cuts \eqref{eq: lbbd_cut}, the BDD-based Benders cuts \eqref{eq: bdd_cut}, and the classical Benders cuts \eqref{eq: bd_cut},  through branch-and-cut. This is different from the standard cutting plane implementation, where master problem need to be solved again each time when new cutting planes are generated. 

The overall implementation of the two-stage decomposition algorithm is illustrated as a flow chart in the Online Supplement. 
In the rest of this section we describe the overall implementation of our decomposition algorithm. We divide the algorithm into two phases for the ease of explanation. In phase one we use the adapted FFD and subproblem relaxation to tighten the UB and LB of the master problem, before entering phase two where we run the branch-and-cut algorithm:

\underline{\it Phase one:} First, we use the adapted FFD heuristic to obtain an initial solution. This solution is added as a {\it warm start} in the commercial solver to provide a feasible solution before the start of the branch-and-cut algorithm. We also generate LBBD cuts or BDD-based Benders cuts, as well as classical Benders cuts from this solution and add them to the master problem. Next, we generated the constraints \eqref{eq: qlb} from the subproblem relaxations and also add them to the master problem.  

\underline{\it Phase two:} We obtain the master problem with extra cuts and constraints from phase one and solve it with branch-and-cut. In the branch-and-cut algorithm at each branch-and-bound node, some integer variables are restricted while the others are linearly relaxed to obtain a node LP relaxation. Solve this LP. If the objective value is greater than or equal to the incumbent UB, then the current node can be pruned. Otherwise we proceed to check the integrality of $(\hat{u}, \hat{y}, \hat{x}, \hat{w})$ in the master solution. If $(\hat{u}, \hat{y}, \hat{x}, \hat{w})$ is integral, then we solve the corresponding subproblems \eqref{model: sors_sub} and subproblem LP relaxations, and get their respective optimal objective values $\bar{Q}$ and $\bar{Q}^{\text{LP}}$. We check if the incumbent UB can be updated with a heuristic solution as described in Section \ref{chapter: heurcb}. Also, we generate LBBD cuts or BDD-based Benders cuts and classical Benders cuts. In the CPLEX lazy constraint callback, compare the master solution of $\hat{Q}$ with $\bar{Q}$ and $\bar{Q}^{\text{LP}}$. If $\hat{Q} < \bar{Q}$, we add LBBD cuts or BDD-based Benders cuts; if $\hat{Q} < \bar{Q}^{\text{LP}}$, we add classical Benders cuts. On the other hand, if some elements in the solution $(\hat{u}, \hat{y}, \hat{x}, \hat{w})$ are fractional, we only solve the subproblem LP relaxations, obtain $\bar{Q}^{\text{LP}}$, generate the classical Benders cuts and implement them if {\cred $\hat{Q} < \bar{Q}^{\text{LP}}$} within the CPLEX user cut callback. We refer to those classical Benders cuts added at fractional solutions as {\it user cuts}. Notice that in order to avoid adding too many user cuts to the point of slowing down the algorithm, we do the following {\it user cut management}: we only add user cuts at the root node and every 150 nodes afterwards, and after 4000 nodes are processed in the branch-and-bound tree, we no longer add any user cut. In addition, we add at most 50 user cuts at the root node and at most 5 user cuts at any node afterwards.

After cutting planes are added in the CPLEX lazy constraint callback or the user cut callback, the node LP relaxation is solved again with those additional cutting planes. We repeat this process, until the stopping criteria is met, i.e. the gap between branch-and-bound UB and LB is small enough. In our implementation we stop the algorithm when such a gap is no more than 1\%.
\section{Three-stage Decomposition}\label{chapter: algorithm2} 
As an alternative to the two-stage decomposition provided in Section \ref{chapter: algorithm1}, we also develop a three-stage decomposition for \DE. Here, we first decompose it into a master problem and a set of subproblems, and develop a class of LBBD optimality cuts to connect them. However, this time each subproblem is more complex, because it is in fact itself a 2SIP. Therefore, we further develop a two-stage decomposition for the subproblems. 

In Section \ref{chapter: lbbd2}, we introduce the LBBD algorithm for the (DE). In Section \ref{chapter: bdd2}, we show how the computationally expensive LBBD subproblems can be further decomposed via the BDD-based method. Additional algorithmic enhancements, some of which {\cred are} similar to those in the two-stage decomposition, are also used in the three-stage decomposition, which we cover in Section \ref{chapter: algo_enhancements2}.

\subsection{Decomposition of the SAA Problem}\label{chapter: lbbd2}
We first decompose (DE) into two stages that are linked via LBBD cuts. In order to distinguish this decomposition with the one in Section \ref{chapter: bdd2}, we use the conventional names {\it LBBD master problem} and {\it LBBD subproblem} for its master and subproblems, respectively.
\subsubsection{LBBD Decomposition Framework.}
The LBBD master problem includes the original variables $u_{hd}$ and $w_p$, the new variables $y_{hd}$, $x_{hdp}$, and a new variable $Q_{hd}$. That is, in the master problem we only decide if patients should be assigned to a hospital on a specific day, leaving the patient to room assignment to the subproblems. Thus, we drop the indices $r\in\mc{R}_h$ in the original variables, and get new variables $y_{hd}$ and $x_{hdp}$. Note that $y_{hd}$ is a nonnegative integer variable, representing the number of {\cred ORs opened in} hospital $h$ on day $d$. Also we use the variable $Q_{hd}$ to denote the {\it expected} cancellation cost for an $(h, d)$ pair:
 \begin{subequations}\label{model: sors2_lbbd_master}
\begin{alignat}{5}
\min~~&\sum_{h\in\mc{H}} \sum_{d\in\mc{D}} G_{hd}u_{hd}+\sum_{h\in\mc{H}}\sum_{d\in\mc{D}}F_{hd}y_{hd}+\sum_{h\in\mc{H}}\sum_{d\in\mc{D}}\sum_{p\in\mc{P}}{\cred c^{\text{sched}}_{dp}}x_{hdp}\hspace{-1cm}\nonumber\\
&+ \sum_{p\in\mc{P} \setminus \mc{P'}}c^{\text{unsched}}_p w_p+ \sum_{h\in\mc{H}}\sum_{d\in\mc{D}} Q_{hd} \label{eq: sors2_master_obj}\\
\st ~~& \sum_{h\in\mc{H}}\sum_{d\in\mc{D}}x_{hdp}=1 & & \forall p \in \mc{P'} \label{eq: sors2_master_assign_mandatory}\\
	  & \sum_{h\in\mc{H}}\sum_{d\in\mc{D}}x_{hdp}+w_p=1 & & \forall p \in \mc{P}\setminus\mc{P'} \label{eq: sors2_master_assign_nonmandatory}\\
     & y_{hd}\leq |\mc{R}_h| u_{hd} & &\forall h\in \mc{H}, d\in\mc{D} \label{eq: sors2_master_y-u}\\
     &  y_{hd}\geq x_{hdp}&& \forall h\in\mc{H}, d\in\mc{D}, p\in\mc{P} \label{eq: sors2_master_y-x}\\
     & [\text{LBBD cuts}] \label{eq: sors2_master_cuts}\\
     & u_{hd},x_{hdp}\in \{0,1\}& &\forall h\in \mc{H}, d\in\mc{D}, p\in \mc{P}, r\in\mc{R}_h\\
     & w_p\in\{0,1\} & &\forall p\in\mc{P}\\
     & y_{hd}\in \Z^+, &&\forall h\in \mc{H}, d\in\mc{D}\\
     &Q_{hd}\geq 0, &&\forall h\in\mc{H}, d\in\mc{D}
\end{alignat}
\end{subequations}

Constraints \eqref{eq: sors2_master_assign_mandatory} and \eqref{eq: sors2_master_assign_nonmandatory} assign patients to hospitals, as well as to days of the current planning horizon or a future time. Constraints \eqref{eq: sors2_master_y-u} ensure that when the surgical suite of hospital $h$ is not opened on a day $d$, there is no OR in that hospital to be opened. Constraints \eqref{eq: sors2_master_y-x} make sure that when a hospital $h$ has a patient assigned on day $d$, there must be at least one OR opened in that hospital. Those constraints are not necessary for a correct formulation, but they make the formulation tighter. Constraints \eqref{eq: sors2_master_cuts} are LBBD cuts that are generated progressively from the LBBD subproblems.

In this framework we have one subproblem per $(h,d)$ pair, minimizing the expected cancellation cost by finding out the best way to assign patients selected by the master problem to the opened ORs. Via solving the master problem, we get the optimal solution ($\hat{u}_{hd}$, $\hat{y}_{hd}$, $\hat{x}_{hdp}$, $\hat{w}_p$, $\hat{Q}_{hd}$). Then we pass $\hat{x}_{hdp}$ and $\hat{y}_{hd}$ to the subproblem. The LBBD subproblem for $(h,d)$ is as below:
 \begin{subequations}\label{model: sors2_sub}
\begin{alignat}{5}
{\cred \mc{Q}_{hd}(\hat{x}_{hd\cdot},\hat{y}_{hd}, T^s_\cdot)}= \min~~& \frac{1}{|\mc{S}|}\sum_{s\in\mc{S}}\sum_{p\in \mc{P}}\sum_{r\in\mc{R}_{h}} c^{\text{cancel}}_p (x_{pr}-z^s_{pr})\hspace*{0cm}&&& \label{eq: sors2_sub_obj}\\
\st ~~& \sum_{r\in\mc{R}_h}x_{pr}=\hat{x}_{hdp} & & \forall p \in \mc{P} \label{eq: sors2_sub_assign_patient}\\
     & \sum_{p\in \mc{P}}T^s_p z^s_{pr} \leq B_{hd} & & \forall r\in\mc{R}_h, s\in\mc{S} \label{eq: sors2_sub_time_limit}\\
     & z^s_{pr} \leq x_{pr} & & \forall p\in \mc{P}, r\in\mc{R}_h, s\in\mc{S} \label{eq: sors2_sub_z-x}\\
     &x_{pr}=0 & &\forall r\geq \hat{y}_{hd}+1 \label{eq: sors2_sub_y-sym_break}\\
     & x_{pr},z^s_{pr}\in \{0,1\}& &\forall p\in \mc{P}, r\in\mc{R}_h, s\in\mc{S}
\end{alignat}
\end{subequations}

Constraints \eqref{eq: sors2_sub_assign_patient} require that the patient to OR assignment for the current $(h,d)$ pair should be consistent with the assignment decision in the master problem. Constraints \eqref{eq: sors2_sub_time_limit} specify time limit. Constraints \eqref{eq: sors2_sub_z-x} ensure that if a patient is not scheduled then she will not be operated on. Constraints \eqref{eq: sors2_sub_y-sym_break} are symmetry breaking constraints for ORs. Those constraints ensure that for a hospital-day pair, ORs with lower indices are opened before ORs with higher indices. These are valid because ORs in the same hospital are homogeneous. Also note that those constraints are connected to the master problem by $\hat{y}_{hd}$, which means at most $\hat{y}_{hd}$ many ORs can be used. Moreover, as the master tries to minimize opened ORs due to opening costs, at the optimal solution, the subproblems will use the exact number of ORs suggested by the master.
\subsubsection{LBBD cuts.}
The LBBD subproblem is always feasible, therefore there is no need to introduce LBBD feasibility cuts. {\cred Using $\bar{Q}_{hd}$ to denote $\mc{Q}_{hd}(\hat{x}_{hd\cdot},\hat{y}_{hd}, T^s_\cdot)$, when} $\hat{Q}_{hd}<\bar{Q}_{hd}$, we add the followings cuts, which we refer to as the LBBD optimality cut, to the LBBD master problem:
 \begin{subequations}\label{eq: lbbd_cut2}
\begin{alignat}{5}
&Q_{hd}\geq \bar{Q}_{hd}\left(g_{hdj}-\sum_{p\in\hat{\mc{P}}_{hd}}(1-x_{hdp})\right) \label{eq: LBBD3_1}\\
&y_{hd}\geq (1+\hat{y}_{hd})(1-g_{hdj}) \label{eq: LBBD3_2}\\
&g_{hdj}\in\{0,1\}
\end{alignat}
\end{subequations}
where $\hat{\mc{P}}_{hd}=\{p\in\mc{P}~|~\hat{x}_{hdp}=1\}$ is the current patient list, index $j=\hat{y}_{hd}$, and $g_{hdj}$ is a binary variable that equals 1 when $y_{hd}\leq \hat{y}_{hd}$, otherwise 0. 

This cut is deducted from the logic $y_{hd}\leq\hat{y}_{hd}=>Q_{hd}\geq \bar{Q}_{hd}-\sum_{p\in\hat{\mc{P}}_{hd}}\bar{Q}_{hd} (1-x_{hdp})$, which means the LB on $Q_{hd}$ is imposed only when $y_{hd}\leq\hat{y}_{hd}$. The lower-bounding inequality ensures that when the current patient list or a superset of it occurs, the corresponding cancellation cost will not be lower than $\bar{Q}_{hd}$. We state the validity of the LBBD cut in Theorem \ref{theorem: lbbd2_valid}, of which the proof is provided in the Online Supplement. 

\begin{theorem}\label{theorem: lbbd2_valid}
The LBBD optimality cuts \eqref{eq: lbbd_cut2} are valid.
\end{theorem}

{\cred A nice feature of cuts \eqref{eq: lbbd_cut2} is that at most $\sum_{h\in\mc{H}}\sum_{d\in\mc{D}}(|\mc{R}_h| + 1)$ many binary variables $g_{hdj}$ are needed, {\it regardless} of the number of LBBD cuts that are eventually generated. Therefore, in our implementation we create all potential variables $\{g_{hdj}|h\in\mc{H}, d\in\mc{D}, j = 0, 1,...,|\mc{R}_h|\}$ before the algorithm starts. This is important because the solver, CPLEX, does not allow the addition of new variables during the branch-and-bound process. This also means we do not need to add too many new variables in order to use these cuts. }
\subsection{Decomposition of the LBBD Subproblems}\label{chapter: bdd2}
The LBBD {\cred subproblem \eqref{model: sors2_sub}} is itself a 2SIP, whose first stage decision is to assign the patients in the list $\hat{\mc{P}}_{hd}$ to ORs, which corresponds to the decision variables $x_{pr}$. The recourse decision concerns whether to operate a surgery under each scenario, represented by the variables $z_{pr}^s$. Solving such a 2SIP can be time consuming. For example in our experiment, for an instance with $|\mc{S}| = 50, |\mc{P}| = 20, |\mc{H}| = 3, |\mc{R}_h| = 3$, when  the LBBD master problem assigns all patients to a single $(h, d)$ pair, the subproblem 2SIP could not be solved to optimality by the commercial solver after one hour. 

Because the subproblem is a 2SIP, we are able to use BDD-based Benders cuts in the decomposition, following a similar approach as in Section \ref{chapter: bdd}. To distinguish from the LBBD decomposition in the previous section, we name the master problem and a subproblem in this decomposition as the {\it BDD master problem} and {\it BDD subproblem}, respectively.

The BDD master problem assigns patients to ORs in an $(h, d)$ pair, and relaxes the constraints on the OR operating time limit. The continuous variabls $\theta_{sr}$ corresponds to the optimal objective value of the subproblem for OR $r$ in scenario $s$. The BDD master problem is as follows: 
\begin{subequations}\label{model: bdd2_master}
\begin{alignat}{5}
\min~~& \sum_{p\in \mc{P}}\sum_{r\in\mc{R}_{h}} c^{\text{cancel}}_p x_{pr}+\frac{1}{|\mc{S}|}\sum_{s\in\mc{S}}\sum_{r\in\mc{R}_{h}} \theta_{sr} \hspace*{0cm}&&&\\
\st ~~& \sum_{r\in\mc{R}_h}x_{pr}=\hat{x}_{hdp} & & \forall p \in \mc{P} \\
     &x_{pr}=0 & &\forall r\geq \hat{y}_{hd}+1\\
     &-\sum_{p\in\mc{P}}c^{\text{cancel}}_p x_{pr} \leq \theta_{sr}\leq 0 &&\forall s\in\mc{S}, r\in \mc{R}_h \label{eq: bdd_master_thetaLB}\\
     & [\text{BDD-based Benders cuts}]\label{eq: bdd_master_cuts}\\
     & x_{pr}\in \{0,1\}& &\forall p\in \mc{P}
\end{alignat}
\end{subequations}
The left-hand side of the equation \eqref{eq: bdd_master_thetaLB} is derived from the relaxation of BDD subproblem \eqref{model: bdd2_sub}, by relaxing its time limit constraint. Constraints \eqref{eq: bdd_master_cuts} are BDD-based Benders cuts generated from the shortest path problem on a BDD transformed from subproblem \eqref{model: bdd2_sub}. All the other constraints have the same meaning as their counterparts in LBBD subproblem \eqref{model: sors2_sub}.

We use $\check{x}_{pr}$ to denote the optimal BDD master solution for $x_{pr}$ and pass it to the BDD subproblem. We observe that the subproblem decomposes with scenarios and ORs, yielding the BDD subproblem for each scenario $s$ and OR $r$ as:
\begin{subequations}\label{model: bdd2_sub}
\begin{alignat}{5}
{\cred \theta_{sr}(\check{x}_{\cdot r}, T^s_\cdot)}= \min~~& \sum_{p\in \mc{P}} - c^{\text{cancel}}_p z^s_{pr} \hspace*{3cm}\\
\st ~~& \sum_{p\in \mc{P}} T^s_p z^s_{pr} \leq B_{hd} && \label{eq: bdd2_sub_timelim} \\
     & z^s_{pr} \leq \check{x}_{pr} & & \forall p\in \mc{P}\label{eq: bdd2_sub_zUB} \\
     & z^s_{pr}\in \{0,1\}& &\forall p\in \mc{P}
\end{alignat}
\end{subequations}
Constraints \eqref{eq: bdd2_sub_timelim} and \eqref{eq: bdd2_sub_zUB} are respectively the time limit constraints and the constraints for UB of $z^s_{pr}$. {\cred Like before, we simplify the notations and use $\ddot{\theta}_{sr}$ to denote $\theta_{sr}(\check{x}_{\cdot r}, T^s_\cdot)$.}

Observe that this problem is almost the same as problem \eqref{model: sors_sub}, except for the constant term in the objective and hospital-day indices. Therefore, we can derive strengthened BDD-based Benders cuts similar to \eqref{eq: bdd_cut}. 

\subsection{Additional Algorithmic Enhancements for the Three-stage Decomposition}\label{chapter: algo_enhancements2}
In this section we introduce several algorithmic enhancements for the three-stage decomposition. In Section \ref{chapter: benders2} we derive classical Benders cuts from LBBD and BDD subproblems. Section \ref{chapter: early_stop} introduces a scheme that avoids solving some LBBD subproblems to optimality by identifying suboptimal master solutions early on. Sections \ref{chapter: ffd2} and \ref{chapter: qlb2} show how to implement an adapted FFD heuristic and add the relaxed subproblem to the master problem in the three-stage decomposition. Section \ref{subch: overall_implement2} explains the overall implementation of the three-stage decomposition algorithm with enhancements. Note that in the three-stage decomposition we also provide additional heuristic solutions in the branch-and-cut implementation the same way as described in Section \ref{chapter: heurcb}, which we will also discuss in Section \ref{subch: overall_implement2}.
\subsubsection{Classical Benders Cuts.}\label{chapter: benders2}
For both the LBBD subproblem and the BDD subproblem, we can derive classical Benders cuts from their LP relaxations, respectively providing LBs for $Q_{hd}$ and $\theta_{sr}$. The procedure of deriving those classical Benders cuts is standard. We describe the process and the classical Benders cuts in the Online Supplement. 
\subsubsection{Early Stopping Scheme in the LBBD Branch-and-cut Implementation.}\label{chapter: early_stop}

We apply an early stopping scheme \citep{karwan1976surrogate} which identifies a suboptimal LBBD master problem solution early on, so we can avoid solving the corresponding 2SIP LBBD subproblems to optimality. This is useful as the LBBD subproblems can be very time consuming to solve, especially considering the LBBD master problem assigns patients to $(h, d)$ pairs without knowledge about the operating time limit of ORs. Such lack of knowledge sometimes leads to the LBBD master problem generating suboptimal solutions that over-schedule patients in some $(h, d)$ pairs and results in high cancellation costs. We would like to identify those suboptimal solutions as soon as possible, because when the LBBD master problem passes a heavily over-scheduled patient list to an $(h, d)$ pair, it can be time-consuming to solve the corresponding subproblem. 

To understand how the early stopping scheme works, first note that we embed the LBBD cuts into the branch-and-bound tree, i.e. in a branch-and-cut framework. While exploring the branch-and-bound tree of the LBBD master problem, we solve an LBBD subproblem at every integral node. For each integral solution, we record the current best UB of the LBBD master problem as {\it globalUB} and the operational cost at the current integral solution as {\it incumbentOptCost}. When solving the LBBD subproblem at such a node, we record the best LB of the subproblem whenever one is found, as $\bar{Q}_{hd}^{\text{LB}}$. If {\it globalUB} $<$ {\it incumbentOptCost}$+ \bar{Q}_{hd}^{\text{LB}}$, there is no need to continue solving the subproblem, as the current LBBD master solution cannot be optimal.

In short, we stop solving an LBBD subproblem once we detect that its lower bound is high enough, and will result in a higher master problem objective value than the best-known UB of the master problem in the branch-and-bound tree. 



Although when the subproblem is stopped early it is not solved to optimality, it is still helpful to derive an LBBD cut from its best-known LB, $\bar{Q}_{hd}^{\text{LB}}$. Such a cut is useful to cut off the current master solution, and has the same logic as \eqref{eq: lbbd_cut2}, though we need to replace $\bar{Q}_{hd}$ in \eqref{eq: lbbd_cut2} with $\bar{Q}_{hd}^{\text{LB}}$:
\begin{subequations}\label{eq: ES_LBBD_cuts}
\begin{alignat}{5}
&Q_{hd}\geq \bar{Q}_{hd}^{\text{LB}}\left(g_{hdj}-\sum_{p\in\hat{\mc{P}}_{hd}}(1-x_{hdp})\right) \\
&y_{hd}\geq (1+\hat{y}_{hd})(1-g_{hdj}) \\
&g_{hdj}\in\{0,1\}
\end{alignat}
\end{subequations}
which tells the LBBD master problem that for the current $(h,d)$ pair, if we open no more than $\hat{y}_{hd}$ rooms and assign all the patients in $\hat{\mc{P}}_{hd}$ to it, the cancellation cost is at least $\bar{Q}_{hd}^{\text{LB}}$. 


\subsubsection{Adapted FFD Heuristic.}\label{chapter: ffd2}
We use the adapted FFD heuristic to get an initial solution for the LBBD master problem, following the same steps as described in Section \ref{chapter: ffd}. We also develop an adapted FFD heuristic for the BDD master problem to assign patients in $\hat{\mc{P}}_{hd}$ to the $\hat{y}_{hd}$ opened ORs. Since the ORs are homogenous in the same $(h, d)$, without loss of generality we select ORs with the smallest $\hat{y}_{hd}$ indices to open, and sort them in the order of their indices. We sort the patients in $\hat{\mc{P}}_{hd}$ by the decreasing order of their cancellation costs, and assign patients to ORs in their sorted order. Based on the surgery duration in the first scenario, for each patient we try to fit her into one of the opened ORs. If none of the opened ORs can fit her, we use an OR that has not been opened. If a patient cannot be fitted into any of the ORs and all the ORs have been used, we assign her to the last OR in the sorted order. After assigning all the patients, we solve a BDD subproblem with this assignment to find the corresponding values of $\ddot{\theta}_{sr}$. This heuristic solution is used before the branch-and-cut process as a warm start. We also derive BBD-based Benders cuts and classical Benders cuts based on the heuristic solution and add them to the BBD master problem \eqref{model: bdd2_master}. 
\subsubsection{Adding Subproblem Relaxations to the Master Problem.}\label{chapter: qlb2}
We derive the following constraints from the LBBD subproblems and add them to the LBBD master problem, to provide lower bounds for the variables $\{Q_{hd}\}_{h\in\mc{H}, d\in\mc{D}}$:
\begin{align}\label{eq: qlb2}
Q_{hd} \geq \frac{1}{|\mc{S}|} \sum_{s\in\mc{S}} \left(\min_{p\in\mc{P}} \frac{c_p^{\text{cancel}}}{T^s_p}\right)\left(\sum_{p\in\mc{P}} T^s_p x_{hdp} - B_{hd} y_{hd}\right) ~~~\forall h\in\mc{H}, d\in\mc{D}
\end{align}
and those constraints are valid:
\begin{theorem}\label{theorem: qlb2}
The constraints \eqref{eq: qlb2} are valid for (DE).
\end{theorem}
The proof is provided in the Online Supplement. 
\subsubsection{Overall Implementation.}\label{subch: overall_implement2}
Similar to the two-stage decomposition, we use a branch-and-cut approach in both the LBBD decomposition and the decomposition of LBBD subproblem. More details on the implementation are provided in the Online Supplement. 

\section{Computational Results}\label{chapter: result}
In this section, we present the computational analysis of the SDORS model and the decomposition algorithms. {\cred Section \ref{subch: parameter} provides the information about our data source and the setup of parameters. Section \ref{subch: saa_analysis} includes the results of the SAA analysis, which decides the number of scenarios to use in (DE). Section \ref{subch: algo_compare} compares the performance of our algorithms with that of the CPLEX solver. Section \ref{subch: incorporate_sto} evaluates the value of incorporating stochasticity in the model. Section \ref{subch: sens_analysis} presents results from the sensitivity analysis.}
\subsection{Parameter Setup}\label{subch: parameter}
In our analysis we use the {\cred same dataset as \cite{roshanaei2017propagating}, which is extracted from the information of 7500 surgeries between 2011 and 2013 from General Surgery Departments of UHN. More specifically, we generate the same parameters $B_{hd}$, $F_{hd}$, $G_{hd}$, {\cred $c^{\text{sched}}_{dp}$ and $c^{\text{unsched}}_{dp}$, which are presented in the 
Online Supplement. The only difference in our setup} lies in the surgery durations and the cost {\cred parameter} $c^{\text{cancel}}_p$, as described below.}

We assume all surgeries are identical, and surgery durations follow a truncated lognormal distribution \citep{strum2000modeling} with mean equals 160 minutes, standard deviation of 40 minutes, and the durations are truncated at 45 minutes and 480 minutes. Those parameters for the surgery duration distribution come from a surgery dataset of UHN.

{\cred The health status score of patient $p$, $\omega_p$, is calculated by $(\alpha_p - |\mc{D}|)\rho_p$, where $\alpha_p$ is the number of days the patient has been waiting since her referral date, and $\rho_p$ is the patient's urgency score. }The benefit of scheduled patients, $c^{\text{sched}}_{dp}$, the penalty for unscheduled patients, $c^{\text{unsched}}_{dp}$, and the penalty for cancelled patient, $c^{\text{cancel}}_p$, are calculated as follows: let $\kappa_1$ be the unit benefit of scheduling patients on time, $\kappa_2$ be the unit penalty of postponing patients to future planning horizons, while $\kappa_3$ and $\kappa_4$ be the unit penalty of canceling non-mandatory and mandatory patients, respectively. Note that we have $\kappa_1>0>\kappa_2>>\kappa_3>\kappa_4$. Then $c^{\text{sched}}_{dp} = \kappa_1\rho_p(d-\alpha_p)$, $c^{\text{unsched}}_{p} = \kappa_2\rho_p(|\mc{D}|+1-\alpha_p)$, $c^{\text{cancel}}_{p} = \kappa_3\rho_p(|\mc{D}|+1-\alpha_p)$ for non-mandatory patients, and $c^{\text{cancel}}_{p} = \kappa_4\rho_p(|\mc{D}|+1-\alpha_p)$ for mandatory patients. {\cred Notice that the benefit of scheduling patients $c^{\text{sched}}_{dp}$ depends on the day of scheduling $d\in\mc{D}$, which means the earlier a patient is scheduled, the larger this benefit is. $c^{\text{sched}}_{dp}\leq 0$ because it represents the benefit while the model's objective is to minimize the total cost. Since $\alpha \in [60, 120]$, a smaller $d$ will lead to a larger $|c^{\text{sched}}_{dp}|$. We use $|\kappa_1| = \$50$ and $|\kappa_2|=\$5$, same as \cite{roshanaei2017propagating}. We cannot directly obtain the cost of cancellation to each patient, so we use the values of $\kappa_1$ and $\kappa_2$ to estimate these penalties. For the non-mandatory patients' unit cancellation costs, we set $|\kappa_3|= \$80 >|\kappa_1|$, so when a surgery is cancelled, the unit loss will be higher than the benefit of scheduling it. The unit cost of cancelling a mandatory patient is $|\kappa_4|= \$100 > |\kappa_3|$. In Section \ref{subch: sens_analysis} we conduct sensitivity analysis on the  cancellation cost.}

{\cred The details for the values of parameters} are provided in the Online Supplement. 

\subsection{The SAA Analysis}\label{subch: saa_analysis}
To determine the number of scenarios to use in the (DE), we derive statistically valid LBs and UBs on the optimal value of the original stochastic program, and compare the worst case optimality gap.

Let us rewrite the original 2SIP in a more concise way:
\begin{align}\label{model: 2sp}
\min_{X\in \mc{X}} ~f(X)=c^\top X+\Exp_{\xi}[Q(X,\xi)]
\end{align}	
where the set $\mc{X}$ include all feasible first-stage variables $X = \{u_{hd}, y_{hdr}, x_{hdpr}, w_p\}$, $c$ represents their objective coefficients, $\xi$ is the vector of random variables that follows a truncated normal distribution, and $Q$ is the second stage value function. 

Considering that we use $|\mc{S}|$ samples in \DE, in order to obtain a LB on the optimal value of problem \eqref{model: 2sp} via SAA, we first generate $N^{\text Sample}$ samples where each sample contains $|\mc{S}|$ scenarios of simulated surgery durations. Due to the limit of computational resource, $N^{\text Sample}$ can be relatively small. For each sample $n$, we solve the (DE) and obtain the optimal solution and objective value. We denote the optimal first-stage solution by $X_n$ and save \sout{\cred the} it because it will be useful when calculating the UB. We denote the optimal objective value by $F^{\text {LB}}_n$. We can then construct a 95\% confidence interval (CI) for the LB of \eqref{model: 2sp} from the $N^{\text{Sample}}$ optimal objective values, $\{F^{\text{LB}}_n\}_{n = 1,...,N^{\text{Sample}}}$. 

To get the UB of SAA, we start with picking a feasible first-stage solution $X\in \mc{X}$. The goal is to choose a solution that has the potential to give a better (i.e. lower) objective value. A heuristic for selecting the solution is to choose from the solutions $X_n, n = 1,..., N^{\text{Sample}}$ as follows: First, generate a medium-sized set of scenarios $\mc{S}^{\text{Select}}$. Then for each solution $X_n$, calculate the second-stage cancelation cost under each scenario ${\cred s\in\mc{S}^{\text{Select}}}$, denoted by $Q(X_n, \xi_s)$. Choose the $X_n$ with the smallest value of $c^\top X_n + \frac{1}{|\mc{S}^{\text{Select}}|} \sum_{s\in\mc{S}^{\text{Select}}} Q(X_n, \xi_s)$ and denote it by $X^{min}$.
 
 After fixing the first-stage solution to $X^{\min}$, we generate a set of scenarios, $\mc{S}^{\text{UB}}$. The cardinality of this set should be large. Then for each scenario $s$ we evaluate the cancellation cost and denote it by $Q(X^{\text{min}}, \xi_s)$. Since $X^{\text{min}}\in\mc{X}$, the optimal objective value from each scenario, $c^\top X^{\text{min}} + Q(X^{\text{min}}, \xi_s)$, is an UB for the optimal objective value of the corresponding scenario problem. We build a 95\% CI from those UBs.

We choose one instance with $|\mc{P}| = 12, |\mc{H}| = 2, |\mc{D}| = 2, |\mc{R}_h| = 2$ to conduct the SAA analysis. As for the parameters of the SAA analysis we use $N^{\text{Sample}} = 30$, $|\mc{S}^{\text{Select}}| = 1,000$ and $|\mc{S}^{\text{UB}}|$ = 10,000. We compare the worst case optimality gap, obtained as ((Mean of UB + Width of UB) - (Mean of LB - Width of LB)) / (Mean of LB - Width of LB), when $|\mc{S}| = 25, 50, 75, 100$. The results are shown in Table \ref{table: saa_results}. 

\begin{table}[h]
\small
\centering
\caption{SAA optimality gap for different levels of scenarios}
\label{table: saa_results}
\begin{tabular}{rrrrrr}
\hline
\multicolumn{1}{c}{}                          & \multicolumn{2}{c}{95\% CI on LB}                    & \multicolumn{2}{c}{95\% CI on UB}                    & \multicolumn{1}{c}{Worst Case}  \\ \cline{2-3}\cline{4-5}
\multicolumn{1}{c}{$|\mc{S}|$} & \multicolumn{1}{c}{Mean} & \multicolumn{1}{c}{Width} & \multicolumn{1}{c}{Mean} & \multicolumn{1}{c}{Width} & \multicolumn{1}{c}{Opt Gap(\%)} \\ \hline
25  & -127363      & 239.67 & -124145 & 125.00 & 2.81 \\
50  & -126633    & 151.93 & -124577 & 109.22 & 1.83 \\
75  & -126175      & 131.77 & -124405 & 112.17  & 1.59 \\
100 & -125967 & 133.24 & -124510 & 83.22 & 1.33\\\hline
\end{tabular}
\end{table}

From the results, we can see that when there are 100 scenarios used in the (DE), the optimality gap is below 1.5\%. Considering the extensive computational efforts to solve the \DE, we prefer to use a small set of scenarios that can still produce reasonably low optimality gaps. Therefore, we fix the number of scenarios $|\mc{S}| = 100$ in the experiments for the following sections.
%
%
\subsection{Algorithm Comparison}\label{subch: algo_compare}
In this section we compare the computational time of the CPLEX solver and our decomposition algorithms. For all the experiments we use C++ API for CPLEX 12.8. We use the default settings of the solver, except that the number of threads is set to one in order to use callbacks for the branch-and-cut algorithm. {\cred Note that for the fairness of comparison, we also use one thread when solving directly with CPLEX.} The computer used to perform those experiments runs MacOS and has a {\cred 2.3 GHz} Intel Core i5 processor and a 16GB RAM. The time limit for these experiments is {\cred 3 hours} and the relative MIP gap is set as 1\%.

We generate instances with various sizes: 10 to 75 patients are scheduled to either 2 or 3 hospitals, 3 or 5 days, and 3 or 5 ORs per hospital\footnote{{\ccred Those instances, along with the obtained bounds on the optimal objective values, are available for download at \url{https://sites.google.com/site/mervebodr/home/SDORS_Instances.zip?attredirects=0&d=1}.}}. The selection for the number of hospitals and days are based on practical reasons, as there are 3 hospitals in the UHN, and the planning horizon of one week contains five business days. We have at most 75 patients and 5 ORs per hospital because further increasing their sizes will result in very large optimality gaps at the end of the time limit. That being said, we believe the instances we sample are diversified enough to cover both sparse and dense cases, and are able to provide an overview for the algorithmic performance.

\begin{table}[h]
\centering
\small
\caption{Comparison of Algorithms}
\label{table: algo_compare}
{\cred
\begin{tabular}{crrrrrrrrr}
\hline
instance  & \multicolumn{4}{c}{Time/Gap}                                                                                  & \multicolumn{1}{c}{} & \multicolumn{4}{c}{Number of Nodes}                                                                           \\ \cline{1-5} \cline{7-10} 
(p-h-d-r) & \multicolumn{1}{c}{MIP} & \multicolumn{1}{c}{2-BDD} & \multicolumn{1}{c}{2-LBBD} & \multicolumn{1}{c}{3-LBBD} & \multicolumn{1}{c}{} & \multicolumn{1}{c}{MIP} & \multicolumn{1}{c}{2-BDD} & \multicolumn{1}{c}{2-LBBD} & \multicolumn{1}{c}{3-LBBD} \\ \hline
10-2-3-3  & 2.8\%                   & 1.1\%                     & 2.3\%                      & \textbf{17.08(min)}        &                      & 90445                   & 242846                    & 163242                     & 51307                      \\
25-2-3-3  & 8.7\%                   & 7.7\%                     & \textbf{6.9\%}             & 8.8\%                      &                      & 16904                   & 27100                     & 21800                      & 2400                       \\
10-3-5-3  & 4.5\%                   & 2.6\%                     & 2.9\%                      & \textbf{43.68(min)}        &                      & 47072                   & 177400                    & 160008                     & 276312                     \\
25-3-5-3  & 24.1\%                  & 15.5\%                    & 16.4\%                     & \textbf{9.4\%}             &                      & 30400                   & 23432                     & 27460                      & 7623                       \\
50-3-5-3  & 34.5\%                  & \textbf{14.8\%}           & 29.0\%                     & 48.0\%                     &                      & 3619                    & 6802                      & 5763                       & 0                          \\
75-3-5-3  & 72.1\%                  & \textbf{15.5\%}           & 20.2\%                     & 46.7\%                     &                      & 0                       & 3325                      & 3080                       & 0                          \\
10-2-3-5  & 3.1\%                   & 62.06(min)                & 2.0\%                      & \textbf{15.37(min)}        &                      & 34473                   & 101474                    & 144058                     & 15337                      \\
25-2-3-5  & 7.5\%                   & 7.8\%                     & \textbf{7.2\%}             & -                          &                      & 12800                   & 24417                     & 23366                      & 0                          \\
50-2-3-5  & 61.1\%                  & 14.8\%                    & \textbf{10.6\%}            & -                          &                      & 7467                    & 9212                      & 10023                      & 0                          \\
75-2-3-5  & 53.3\%                  & \textbf{17.3\%}           & 17.5\%                     & -                          &                      & 1                       & 2605                      & 2323                       & 0                          \\
10-3-5-5  & 3.2\%                   & 3.0\%                     & 4.6\%                      & \textbf{18.06(min)}        &                      & 27090                   & 124901                    & 108930                     & 24361                      \\
25-3-5-5  & 28.5\%                  & \textbf{14.4\%}           & 15.5\%                     & -                          &                      & 5000                    & 17894                     & 25506                      & 0                          \\
50-3-5-5  & 59.4\%                  & \textbf{18.1\%}           & 20.6\%                     & -                          &                      & 293                     & 7854                      & 8415                       & 0                          \\
75-3-5-5  & 52.7\%                  & 22.7\%                    & \textbf{16.3\%}            & -                          &                      & 9                       & 3506                      & 3827                       & 0                          \\ \hline
\end{tabular}}
\end{table}

The experimental results are shown in 
Table \ref{table: algo_compare}. The instances are denoted by the number of patients (p), hospitals (h), days (d), and rooms (r), linked by dashes. We solve the \DE~directly with CPLEX (``MIP"), as well as solving it using the two-stage decomposition with either only BDD-based cuts (``2-BDD") or only LBBD cuts (``2-LBBD"), and solving it with the three-stage decomposition (``3-LBBD"). We report either the solution time or the optimality gap at the end of the time limit, along with the number of nodes processed in the branch-and-bound tree. In the columns for ``Time/Gap" comparison, the bold items represent either the smallest solution time or optimality gap, and the ``-" signs represent the cases with over 100\% gap. The results in the ``Number of Nodes" sections show the number of branch-and-bound nodes processed before the algorithms stop. 

The columns of Time/Gap comparisons in Table \ref{table: algo_compare} show that our decomposition algorithms are more efficient than the MIP solver in all test instances. In particular, for the instances with more patients where directly using the MIP solver results in gaps beyond {\cred 50\%,} our decomposition algorithms can solve those instances to a gap {\cred below} 20\%. The three-stage decomposition outperforms the other methods when there are fewer patients to be scheduled. In fact, it performs so well in small instances that in the experiments of Section \ref{subch: saa_analysis} we used this method to solve all the optimization problems to optimality, since the SAA analysis uses a small instance with $|\mc{P}| = 12$. For instances with 25 or more patients, the two-stage decomposition methods are the best. {\cred Between the two two-stage decompositions, one outperforms the other in almost the same number of tested instances. } For the larger instances, the three-stage decomposition becomes very ineffective, as the size of the LBBD subproblems become the bottleneck. 


{\cred We also note that for the two-stage decomposition, it is also possible to use both BDD-based cuts and LBBD cuts at the same time. However, in most of the instances, adding both cuts does not perform as well as the best performing algorithm that has only one of those two types of cuts. As a potential reason of 2-BDD and 2-LBBD collectively outperforming adding both cuts, we observe that adding extra cuts in each iteration quickly increases the number of cuts in the branch-and-cut algorithm, which slows down the algorithm in later rounds, and overrides the benefit of having an extra type of cut. We do acknowledge though that a careful cut management strategy can potentially make adding both cuts a better option.} 

{\cred In addition, we present results for the objective values of the best integer solutions in the Online Supplement. 
The results show that if an algorithm has the best computational performance, it is also most likely to produce the best integer solutions.}
\subsection{Value of Incorporating Stochasticity}\label{subch: incorporate_sto}
We evaluate the value of incorporating stochastic surgery durations in the distributed OR scheduling problem by comparing the optimal schedule from the deterministic {\cred DORS model  \citep{roshanaei2017propagating}} and the ones from our stochastic problem \eqref{model: sdors}. In this experiment, we use those instances with three hospitals and five days from Section \ref{subch: algo_compare}, as those instances are the more realistic ones in terms of the number of hospitals and days. 

For the evaluation, we first solve both models with a time limit of {\cred 3 hours} and the optimality gap being set to 1\%. Note that after {\cred 3 hours}, all the deterministic models are solved to optimality, while for the stochastic models we use the results from the fastest algorithms
as indicated by Table \ref{table: algo_compare}, and most instances are not solved to optimality. After we obtain the hospital and OR opening decisions and patient assignments from the optimization, we evaluate those solutions with randomly generated surgery durations. We randomly generate 10,000 samples of durations, solve the recourse problem \eqref{model: sdors-sec} with the first stage decisions fixed to the values from the optimization, and calculate the {\it cancellation rate} and {\it OR utilization rate} for each sample under different assignments. The cancellation rate equals the ratio between the number of cancelled patients and the total number of scheduled patients. The OR utilization rate is the ratio between the total duration of accepted surgeries and the total available time of all ORs. Finally we {\cred obtain the 95\% confidence intervals for} the evaluations from those 10,000 samples.

\begin{table}[h]
\cred
\small
\centering
\caption{Comparison of the Deterministic and Stochastic Models}
\label{table: deter_sto_compare}
\begin{tabular}{crrrrr}
\hline
instance  & \multicolumn{2}{c}{Cancellation Rate}                                       & \multicolumn{1}{c}{} & \multicolumn{2}{c}{Utilization Rate}                                        \\ \cline{1-3} \cline{5-6} 
(p-h-d-r) & \multicolumn{1}{c}{Deterministic} & \multicolumn{1}{c}{Stochastic}          & \multicolumn{1}{c}{} & \multicolumn{1}{c}{Deterministic} & \multicolumn{1}{c}{Stochastic}          \\ \hline
10-3-5-3  & 18.1\%$\pm$1.9E-3    & \textbf{0.6\%}$\pm$4.9E-4  &                      & 5.6\%$\pm$8.8E-5     & \textbf{7.7\%}$\pm$1.1E-4  \\
25-3-5-3  & 18.9\%$\pm$1.1E-3    & \textbf{14.3\%}$\pm$1.0E-3 &                      & 15.3\%$\pm$1.4E-4    & \textbf{16.2\%}$\pm$1.5E-4 \\
50-3-5-3  & 16.4\%$\pm$7.1E-4    & \textbf{9.6\%}$\pm$6.2E-4 &                      & 31.0\%$\pm$2.1E-4    & \textbf{34.5\%}$\pm$2.3E-4 \\
75-3-5-3  & 17.8\%$\pm$5.9E-4    & \textbf{9.4\%}$\pm$4.6E-4  &                      & 45.9\%$\pm$2.5E-4    & \textbf{52.0\%}$\pm$2.8E-4 \\
10-3-5-5  & 18.2\%$\pm$1.9E-3    & \textbf{0.6\%}$\pm$4.9E-4  &                      & 3.3\%$\pm$5.3E-5     & \textbf{4.6\%}$\pm$6.8E-5  \\
25-3-5-5  & 20.9\%$\pm$1.0E-3    & \textbf{2.5\%}$\pm$5.1E-4  &                      & 8.9\%$\pm$8.3E-5     & \textbf{11.3\%}$\pm$1.1E-4 \\
50-3-5-5  & 16.5\%$\pm$7.2E-4    & \textbf{0.5\%}$\pm$1.9E-4  &                      & 18.5\%$\pm$1.3E-4    & \textbf{23.1\%}$\pm$1.6E-4 \\
75-3-5-5  & 15.3\%$\pm$5.7E-4    & \textbf{4.3\%}$\pm$3.5E-4  &                      & 28.9\%$\pm$1.6E-4    & \textbf{33.2\%}$\pm$1.8E-4 \\ \hline
\end{tabular}
\end{table}

In Table \ref{table: deter_sto_compare} 
we report {\cred 95\% confidence intervals of} the cancellation rate and the OR utilization rate. For each instance, the model with lower cancellation rate and higher OR utilization rate is highlighted. From the table we can conclude that using stochastic model reduces the patient cancellation rate and improves the utilization of ORs. In particular, the improvement is significant for instances with more patients and ORs, e.g., the instances 50-3-5-5 and 75-3-5-5. {\cred This can be explained by the fact that the deterministic model opens an average of 28\% fewer rooms than the stochastic model, while scheduling about the same number of patients.} 

{\cred It is worth mentioning that the utilization rate we report is  calculated considering all available ORs, whether open or not. If we look at the utilization rate of only opened ORs, then {\ccred for the deterministic model the opened OR utilization rates are in the range of 78\% to 83\% for the instances in Table \ref{table: deter_sto_compare}, while for the stochastic model the opened OR utilization rates are between 67\% and 79\%. We observe that compared with its deterministic counterpart, the stochastic model has a lower opened OR utilization rate in all the tested instances, because it opens more ORs in order to reduce the cancellation rate.}

Notice that all the deterministic models are solved to optimality, while most of the stochastic models have an optimality gap of around {\cred 16\%} when the time limit is reached. Nevertheless, the result from a stochastic model provides a more robust schedule. Therefore, it is worth the effort to implement the stochastic model, even if its implementation is more involved.
\subsection{Sensitivity Analysis}\label{subch: sens_analysis}
In this section we perform several sensitivity analysis for the SDORS problem. Specifically, we change the standard deviation of the surgery duration distribution, the cancellation cost, and the operating time limit of ORs, to see how changes in those parameters influence the performance of our model. 

In this experiment it is important to solve all the optimization problems to optimality, because if not, it will be hard to conclude whether the differences in performance are resulted from the change of parameters or differences in optimality gaps. Therefore, all experiments in this section are conducted on {\cred instances} with $|\mc{P}| = 12, |\mc{H}| = 2, |\mc{D}| = 2, |\mc{R}_h| = 2$. {\cred Those instances} can be solved to optimality within a reasonable time (using the three-stage decomposition), but {\cred are} also non-trivial enough to produce interesting results, {\cred where the change of parameters leads to a change in the measures, such as the cancellation rate and the utilization rate. We use ten instances for each of the different parameter settings.}

We first solve the SDORS {\cred instances} with the original parameter setting (``Baseline"). Then we solve the problem again after one of the following changes:

Case 1: Increase the standard deviation of surgery durations from 40 minutes to 60 minutes.

Case 2: Reduce the unit penalty of canceling patients, $\kappa_3$ and $\kappa_4$, to $\frac{2}{3}$ of their original values.

Case 3: Reduce the operating time limits of ORs to half of their original values.

After solving those optimization problems and obtain the optimal schedules, we follow the same procedure as in Section \ref{subch: incorporate_sto} and obtain {\cred the cancellation rate (CancelRate) and utilization rate (UtiliRate) corresponding to each instance. We report the averages and standard deviations of those measure. We also obtain the average numbers and standard deviations of scheduled patients (Scheduled) and opened ORs (OpenOR).} Results are shown in 
Table \ref{table: sens_analysis}{\cred , where the averages of each measure are reported, followed by the standard deviation (STD) in parentheses.}
We compare the baseline case where none of the parameters are changed, {\cred with} the three cases that correspond to the three types of changes in the parameters. 

\begin{table}[h]
\centering
\small
{\cred
\caption{Sensitivity Analysis for SDORS}
\label{table: sens_analysis}
\begin{tabular}{crlrlrlrl}
\hline
Settings & CancelRate & (STD)    & UtiliRate & (STD)    & Scheduled & (STD) & OpenOR & (STD)    \\ \hline
Baseline           & 8.1\%             & (2.9E-4) & 12.6\%           & (4.7E-5) & 12        & (0)   & 5      & (0)      \\
Case 1             & 8.7\%             & (1.9E-2) & 12.3\%           & (3.1E-3) & 12        & (0)   & 5      & (3.2E-1) \\
Case 2             & 16.7\%            & (4.7E-4) & 11.3\%           & (5.6E-5) & 12        & (0)   & 4      & (0)      \\
Case 3             & 5.1\%             & (2.9E-4) & 17.1\%           & (7.5E-5) & 8         & (0)   & 8      & (0)      \\ \hline
\end{tabular}}
\end{table}


%

{\cred From the results, we see that for Case 1 where the standard deviation of the surgery duration becomes higher, we have a higher cancellation rate and a lower OR utilization rate. This is expected, as more uncertainty usually leads to more cancellations. Also, in our experiment, the average number of scheduled patients and the average number of opened ORs are the same for Baseline and Case 1, so a higher cancellation rate naturally leads to a lower OR utilization.

For Case 2, the cancellation rate increases significantly, and the OR utilization rate drops. When cancelling a patient becomes less costly, the hospitals are motivated to schedule patients to fewer ORs to save operational costs. In the experiment, the average number of opened ORs drops from 5 (Baseline) to 4, while the average number of scheduled patients remains the same. This trend explains the increase in the cancellation rate. Also, compared with the drop in opened ORs, the increase in cancellation rate is more drastic, and this leads to an overall effect of a decreased utilization rate. 

For Case 3, the cancellation rate decreases and the OR utilization rate increases. When the opening hours of ORs become shorter, the average number of patients per OR becomes lower. This understandably reduces cancellation, as fewer patients per OR means lower uncertainty. The decrease in the cancellation rate may explain the increase in the OR utilization rate.}

\section{Conclusions}\label{chapter: conclusion}
In this paper we propose the SDORS problem. Through computational experiments we show that, when compared with its deterministic counterpart from the literature, the SDORS model generates more robust schedules that reduce the rate of cancellation and improve the rate of OR utilization.

We also develop several decomposition algorithms and algorithmic enhancements to improve the computational efficiency, and conduct numerical experiments on real instances from the UHN. For all the instances we experimented with, our algorithms are able to reduce the solution time or the optimality gap compared with the commercial solver. Moreover, our algorithms have the potential to be applied to any stochastic distributed bin packing problem, as the SDORS problem can be viewed as an extension of a stochastic distributed bin packing problem.

We recognize that the SDORS problem is a very difficult problem to solve. Although our algorithms perform much better than the commercial solver, there is still room for improvement, to be able to efficiently handle practical instances, where there are hundreds of patients to be scheduled each week. Therefore, future works that further improve the computational efficiency will be valuable. 



\setcitestyle{numbers}
\bibliography{ms} 

\begin{thebibliography}{46}
\providecommand{\natexlab}[1]{#1}
\providecommand{\url}[1]{\texttt{#1}}
\expandafter\ifx\csname urlstyle\endcsname\relax
  \providecommand{\doi}[1]{doi: #1}\else
  \providecommand{\doi}{doi: \begingroup \urlstyle{rm}\Url}\fi

\bibitem[Angulo et~al.(2016)Angulo, Ahmed, and Dey]{angulo2016improving}
G.~Angulo, S.~Ahmed, and S.~S. Dey.
\newblock Improving the integer {L}-shaped method.
\newblock \emph{INFORMS J.~Comput.}, 28\penalty0 (3):\penalty0 483--499, 2016.

\bibitem[Balas et~al.(1993)Balas, Ceria, and Cornu{\'e}jols]{balas1993lift}
E.~Balas, S.~Ceria, and G.~Cornu{\'e}jols.
\newblock A lift-and-project cutting plane algorithm for mixed 0--1 programs.
\newblock \emph{Math.~Program.}, 58\penalty0 (1-3):\penalty0 295--324, 1993.

\bibitem[Barua and Jacques(2019)]{esmail2019private}
B.~Barua and D.~Jacques.
\newblock The private cost of public queues for medically necessary care.
\newblock Technical report, Fraser Institute, 2019.

\bibitem[Behnamian and Ghomi(2013)]{behnamian2013heterogeneous}
J.~Behnamian and S.~F. Ghomi.
\newblock The heterogeneous multi-factory production network scheduling with
  adaptive communication policy and parallel machine.
\newblock \emph{Inf.~Sci.}, 219:\penalty0 181--196, 2013.

\bibitem[Beier et~al.(2015)Beier, Venkatachalam, Corolli, and
  Ntaimo]{beier2015stage}
E.~Beier, S.~Venkatachalam, L.~Corolli, and L.~Ntaimo.
\newblock Stage-and scenario-wise fenchel decomposition for stochastic mixed
  0-1 programs with special structure.
\newblock \emph{Comput. and Oper.~Res.}, 59:\penalty0 94--103, 2015.

\bibitem[Beli{\"e}n and Demeulemeester(2008)]{belien2008branch}
J.~Beli{\"e}n and E.~Demeulemeester.
\newblock A branch-and-price approach for integrating nurse and surgery
  scheduling.
\newblock \emph{European~J.~Oper.~Res.}, 189\penalty0 (3):\penalty0 652--668,
  2008.

\bibitem[Benders(1962)]{benders1962partitioning}
J.~F. Benders.
\newblock Partitioning procedures for solving mixed-variables programming
  problems.
\newblock \emph{Numerische Mathematik}, 4\penalty0 (1):\penalty0 238--252,
  1962.

\bibitem[Bergman et~al.(2016)Bergman, Cire, Van~Hoeve, and
  Hooker]{bergman2016decision}
D.~Bergman, A.~A. Cire, W.-J. Van~Hoeve, and J.~Hooker.
\newblock \emph{Decision diagrams for optimization}.
\newblock Springer, 2016.

\bibitem[Blake and Carter(2002)]{blake2002goal}
J.~T. Blake and M.~W. Carter.
\newblock A goal programming approach to strategic resource allocation in acute
  care hospitals.
\newblock \emph{European~J.~Oper.~Res.}, 140\penalty0 (3):\penalty0 541--561,
  2002.

\bibitem[Blake and Donald(2002)]{blake2002mount}
J.~T. Blake and J.~Donald.
\newblock Mount {S}inai hospital uses integer programming to allocate operating
  room time.
\newblock \emph{Interfaces}, 32\penalty0 (2):\penalty0 63--73, 2002.

\bibitem[Bodur(2015)]{bodur2015valid}
M.~Bodur.
\newblock \emph{On valid inequalities for polyhedra in extended and projected
  spaces with application to two-stage stochastic integer programming}.
\newblock PhD thesis, The University of Wisconsin-Madison, 2015.

\bibitem[Bodur and Luedtke(2016)]{bodur2016mixed}
M.~Bodur and J.~R. Luedtke.
\newblock Mixed-integer rounding enhanced {B}enders decomposition for
  multiclass service-system staffing and scheduling with arrival rate
  uncertainty.
\newblock \emph{Management~Sci.}, 63\penalty0 (7):\penalty0 2073--2091, 2016.

\bibitem[Bowers and Mould(2004)]{bowers2004managing}
J.~Bowers and G.~Mould.
\newblock Managing uncertainty in orthopaedic trauma theatres.
\newblock \emph{European~J.~Oper.~Res.}, 154\penalty0 (3):\penalty0 599--608,
  2004.

\bibitem[Cardoen et~al.(2010)Cardoen, Demeulemeester, and
  Beli{\"e}n]{cardoen2010operating}
B.~Cardoen, E.~Demeulemeester, and J.~Beli{\"e}n.
\newblock Operating room planning and scheduling: A literature review.
\newblock \emph{European~J.~Oper.~Res.}, 201\penalty0 (3):\penalty0 921--932,
  2010.

\bibitem[Car{\o}e and Tind(1998)]{caroe1998shaped}
C.~C. Car{\o}e and J.~Tind.
\newblock L-shaped decomposition of two-stage stochastic programs with integer
  recourse.
\newblock \emph{Math.~Program.}, 83\penalty0 (1-3):\penalty0 451--464, 1998.

\bibitem[CIHI(2020)]{canadian2020wait}
CIHI.
\newblock {Wait times for priority procedures in Canada}.
\newblock
  \url{https://www.cihi.ca/en/wait-times-for-priority-procedures-in-canada},
  2020.
\newblock [Online; accessed October 2nd, 2020].

\bibitem[Cir{\'e} et~al.(2016)Cir{\'e}, Coban, and Hooker]{cire2016logic}
A.~A. Cir{\'e}, E.~Coban, and J.~N. Hooker.
\newblock Logic-based {B}enders decomposition for planning and scheduling: a
  computational analysis.
\newblock \emph{The Knowledge Engineering Review}, 31\penalty0 (5):\penalty0
  440--451, 2016.

\bibitem[Deng et~al.(2019)Deng, Shen, and Denton]{deng2019chance}
Y.~Deng, S.~Shen, and B.~Denton.
\newblock Chance-constrained surgery planning under conditions of limited and
  ambiguous data.
\newblock \emph{INFORMS J.~Comput.}, 31\penalty0 (3):\penalty0 559--575, 2019.

\bibitem[Denton et~al.(2010)Denton, Miller, Balasubramanian, and
  Huschka]{denton2010optimal}
B.~T. Denton, A.~J. Miller, H.~J. Balasubramanian, and T.~R. Huschka.
\newblock Optimal allocation of surgery blocks to operating rooms under
  uncertainty.
\newblock \emph{Oper.~Res.}, 58\penalty0 (4-part-1):\penalty0 802--816, 2010.

\bibitem[Dimitriadis et~al.(2013)Dimitriadis, Iyer, and
  Evgeniou]{dimitriadis2013challenge}
P.~Dimitriadis, S.~Iyer, and E.~Evgeniou.
\newblock The challenge of cancellations on the day of surgery.
\newblock \emph{International~J.~Surg.}, 11\penalty0 (10):\penalty0 1126--1130,
  2013.

\bibitem[Fazel-Zarandi et~al.(2013)Fazel-Zarandi, Berman, and
  Beck]{fazel2013solving}
M.~M. Fazel-Zarandi, O.~Berman, and J.~C. Beck.
\newblock Solving a stochastic facility location/fleet management problem with
  logic-based {B}enders' decomposition.
\newblock \emph{IIE Trans.}, 45\penalty0 (8):\penalty0 896--911, 2013.

\bibitem[Fei et~al.(2008)Fei, Chu, Meskens, and Artiba]{fei2008solving}
H.~Fei, C.~Chu, N.~Meskens, and A.~Artiba.
\newblock Solving surgical cases assignment problem by a branch-and-price
  approach.
\newblock \emph{International~J.~Production~Economics}, 112\penalty0
  (1):\penalty0 96--108, 2008.

\bibitem[Fei et~al.(2009)Fei, Chu, and Meskens]{fei2009solving}
H.~Fei, C.~Chu, and N.~Meskens.
\newblock Solving a tactical operating room planning problem by a
  column-generation-based heuristic procedure with four criteria.
\newblock \emph{European~J.~Oper.~Res.}, 166\penalty0 (1):\penalty0 91, 2009.

\bibitem[Gade et~al.(2014)Gade, K{\"u}{\c{c}}{\"u}kyavuz, and
  Sen]{gade2014decomposition}
D.~Gade, S.~K{\"u}{\c{c}}{\"u}kyavuz, and S.~Sen.
\newblock Decomposition algorithms with parametric gomory cuts for two-stage
  stochastic integer programs.
\newblock \emph{Math.~Program.}, 144\penalty0 (1-2):\penalty0 39--64, 2014.

\bibitem[Guerriero and Guido(2011)]{guerriero2011operational}
F.~Guerriero and R.~Guido.
\newblock Operational research in the management of the operating theatre: a
  survey.
\newblock \emph{Health~Care~Manage~Sci.}, 14\penalty0 (1):\penalty0 89--114,
  2011.

\bibitem[Gul et~al.(2015)Gul, Denton, and Fowler]{gul2015progressive}
S.~Gul, B.~T. Denton, and J.~W. Fowler.
\newblock A progressive hedging approach for surgery planning under
  uncertainty.
\newblock \emph{INFORMS J.~Comput.}, 27\penalty0 (4):\penalty0 755--772, 2015.

\bibitem[Harper(2002)]{harper2002framework}
P.~R. Harper.
\newblock A framework for operational modelling of hospital resources.
\newblock \emph{Health~Care~Manage~Sci.}, 5\penalty0 (3):\penalty0 165--173,
  2002.

\bibitem[Hooker and Ottosson(2003)]{hooker2003logic}
J.~N. Hooker and G.~Ottosson.
\newblock Logic-based {B}enders decomposition.
\newblock \emph{Math.~Program.}, 96\penalty0 (1):\penalty0 33--60, 2003.

\bibitem[Jebali et~al.(2006)Jebali, Alouane, and Ladet]{jebali2006operating}
A.~Jebali, A.~B.~H. Alouane, and P.~Ladet.
\newblock Operating rooms scheduling.
\newblock \emph{International~J.~Production~Economics}, 99\penalty0
  (1-2):\penalty0 52--62, 2006.

\bibitem[Johnson et~al.(1974)Johnson, Demers, Ullman, Garey, and
  Graham]{johnson1974worst}
D.~S. Johnson, A.~Demers, J.~D. Ullman, M.~R. Garey, and R.~L. Graham.
\newblock Worst-case performance bounds for simple one-dimensional packing
  algorithms.
\newblock \emph{SIAM J.~Comput.}, 3\penalty0 (4):\penalty0 299--325, 1974.

\bibitem[Karwan(1976)]{karwan1976surrogate}
M.~H. Karwan.
\newblock \emph{Surrogate constraint duality and extensions in integer
  programming}.
\newblock PhD thesis, Georgia Institute of Technology, 1976.

\bibitem[Laporte and Louveaux(1993)]{laporte1993integer}
G.~Laporte and F.~V. Louveaux.
\newblock The integer {L}-shaped method for stochastic integer programs with
  complete recourse.
\newblock \emph{Oper.~Res.~Lett.}, 13\penalty0 (3):\penalty0 133--142, 1993.

\bibitem[Lombardi et~al.(2010)Lombardi, Milano, Ruggiero, and
  Benini]{lombardi2010stochastic}
M.~Lombardi, M.~Milano, M.~Ruggiero, and L.~Benini.
\newblock Stochastic allocation and scheduling for conditional task graphs in
  multi-processor systems-on-chip.
\newblock \emph{J.~Sched.}, 13\penalty0 (4):\penalty0 315--345, 2010.

\bibitem[Lozano and Smith(2018)]{lozano2018binary}
L.~Lozano and J.~C. Smith.
\newblock A binary decision diagram based algorithm for solving a class of
  binary two-stage stochastic programs.
\newblock \emph{Math.~Program.}, pages 1--24, 2018.

\bibitem[Magnussen et~al.(2007)Magnussen, Hagen, and
  Kaarboe]{magnussen2007centralized}
J.~Magnussen, T.~P. Hagen, and O.~M. Kaarboe.
\newblock Centralized or decentralized? a case study of norwegian hospital
  reform.
\newblock \emph{Social science \& medicine}, 64\penalty0 (10):\penalty0
  2129--2137, 2007.

\bibitem[Marques et~al.(2012)Marques, Captivo, and Pato]{marques2012integer}
I.~Marques, M.~E. Captivo, and M.~V. Pato.
\newblock An integer programming approach to elective surgery scheduling.
\newblock \emph{OR~Spectrum}, 34\penalty0 (2):\penalty0 407--427, 2012.

\bibitem[Min and Yih(2010)]{min2010scheduling}
D.~Min and Y.~Yih.
\newblock Scheduling elective surgery under uncertainty and downstream capacity
  constraints.
\newblock \emph{European~J.~Oper.~Res.}, 206\penalty0 (3):\penalty0 642--652,
  2010.

\bibitem[Naderi and Ruiz(2014)]{naderi2014scatter}
B.~Naderi and R.~Ruiz.
\newblock A scatter search algorithm for the distributed permutation flowshop
  scheduling problem.
\newblock \emph{European~J.~Oper.~Res.}, 239\penalty0 (2):\penalty0 323--334,
  2014.

\bibitem[Ntaimo(2013)]{ntaimo2013fenchel}
L.~Ntaimo.
\newblock Fenchel decomposition for stochastic mixed-integer programming.
\newblock \emph{J.~Global.~Opt.}, 55\penalty0 (1):\penalty0 141--163, 2013.

\bibitem[Ntaimo and Tanner(2008)]{ntaimo2008computations}
L.~Ntaimo and M.~W. Tanner.
\newblock Computations with disjunctive cuts for two-stage stochastic mixed 0-1
  integer programs.
\newblock \emph{J.~Global.~Opt.}, 41\penalty0 (3):\penalty0 365--384, 2008.

\bibitem[Roshanaei et~al.(2017{\natexlab{a}})Roshanaei, Luong, Aleman, and
  Urbach]{roshanaei2017propagating}
V.~Roshanaei, C.~Luong, D.~M. Aleman, and D.~Urbach.
\newblock Propagating logic-based {B}enders decomposition approaches for
  distributed operating room scheduling.
\newblock \emph{European~J.~Oper.~Res.}, 257\penalty0 (2):\penalty0 439--455,
  2017{\natexlab{a}}.

\bibitem[Roshanaei et~al.(2017{\natexlab{b}})Roshanaei, Luong, Aleman, and
  Urbach]{roshanaei2017collaborative}
V.~Roshanaei, C.~Luong, D.~M. Aleman, and D.~R. Urbach.
\newblock Collaborative operating room planning and scheduling.
\newblock \emph{INFORMS J.~Comput.}, 29\penalty0 (3):\penalty0 558--580,
  2017{\natexlab{b}}.

\bibitem[Sherali and Fraticelli(2002)]{sherali2002modification}
H.~D. Sherali and B.~M. Fraticelli.
\newblock A modification of benders' decomposition algorithm for discrete
  subproblems: An approach for stochastic programs with integer recourse.
\newblock \emph{J.~Global.~Opt.}, 22\penalty0 (1-4):\penalty0 319--342, 2002.

\bibitem[Strum et~al.(2000)Strum, May, and Vargas]{strum2000modeling}
D.~P. Strum, J.~H. May, and L.~G. Vargas.
\newblock Modeling the uncertainty of surgical procedure timescomparison of
  log-normal and normal models.
\newblock \emph{Anesthesiol.:~J.~Am.~Soc.~Anesthesiol.}, 92\penalty0
  (4):\penalty0 1160--1167, 2000.

\bibitem[Timpe and Kallrath(2000)]{timpe2000optimal}
C.~H. Timpe and J.~Kallrath.
\newblock Optimal planning in large multi-site production networks.
\newblock \emph{European~J.~Oper.~Res.}, 126\penalty0 (2):\penalty0 422--435,
  2000.

\bibitem[Wang et~al.(2016)Wang, Roshanaei, Aleman, and
  Urbach]{wang2016discrete}
S.~Wang, V.~Roshanaei, D.~Aleman, and D.~Urbach.
\newblock A discrete event simulation evaluation of distributed operating room
  scheduling.
\newblock \emph{IIE~Trans.~Healthc.~Syst.~Engineering}, 6\penalty0
  (4):\penalty0 236--245, 2016.

\end{thebibliography}

\appendix

\section{Proofs}\label{chapter: app_proof}
\medskip
\subsection{Proof of Theorem \ref{theorem: lbbd}}\label{subch: app_proof_th1}
{\cred \begin{theorem}
The LBBD optimality cut \eqref{eq: lbbd_cut}, which is defined as
\begin{align}
{\cred Q^s_{hdr}} \geq \bar{Q}^s_{hdr} -\sum_{p \in \hat{\mc{P}}_{hdr}} c_p^{\text{cancel}}(1 - x_{hdpr}), \nonumber
\end{align}
\end{theorem}
is valid.
}
\proof{Proof.}
Note that the cut \eqref{eq: lbbd_cut} is formulated for one OR, $r$, in one hospital $h$ on one day $d$ for scenario $s$, thus the following argument is formed in terms of a single tuple of $(h, d, r, s)$ and holds for each such tuple.

To prove the validity of the cut, we need to ensure two points: the cut eliminates the current master solution, and it does not exclude any globally optimal solution. 

Let us first prove that the current master solution ($\hat{u}_{hd}$, $\hat{y}_{hdr}$, $\hat{x}_{hdpr}$, $\hat{w}_p$, $\hat{Q}^s_{hdr}$) violates the cut, i.e., it will be cut off. Assume towards contradiction that it satisfies the cut. When we substitute the master solution to the cut \eqref{eq: lbbd_cut}, the left-hand side (LHS) of the cut becomes $\hat{Q}^s_{hdr}$, while the right-hand side (RHS) of the cut becomes $\bar{Q}^s_{hdr}$ as all $\hat{x}_{hdpr}$'s in the set $\hat{P}_{hdr}$ have the value 1. This gives us $\hat{Q}^s_{hdr} \geq \bar{Q}^s_{hdr}$. However, we know that $\hat{Q}^s_{hdr} < \bar{Q}^s_{hdr}$ because otherwise the LBBD optimality cut will not be generated. This is a contradiction, therefore the current master solution does not satisfy the cut.

Next, we need to prove that any globally optimal solution, denoted by ($\hat{u}^*_{hd}$, $\hat{y}^*_{hdr}$, $\hat{x}^*_{hdpr}$, $\hat{w}^*_p$, $\bar{Q}^{s*}_{hdr}$), is not excluded by the LBBD optimality cut that is generated from a master solution ($\hat{u}_{hd}$, $\hat{y}_{hdr}$, $\hat{x}_{hdpr}$, $\hat{w}_p$, $\hat{Q}^s_{hdr}$) using $\hat{P}_{hdr}$. Here, $(\hat{u}^*_{hd}, \hat{y}^*_{hdr}, \hat{x}^*_{hdpr}, \hat{w}^*_p)$ is obtained by solving the master problem, and $\bar{Q}^{s*}_{hdr}$ is the corresponding optimal objective value of the subproblem. For a globally optimal solution, the cancellation cost obtained by the master problem,  $\hat{Q}^{s*}_{hdr}$, should match that from the subproblem,  $\bar{Q}^{s*}_{hdr}$.

We discuss the different cases of globally optimal solutions below:

\underline{\it Case 1}: If $\hat{y}^*_{hdr}=0$, i.e. the OR $r$ in hospital $h$ on day $d$ is closed, then the optimal solution must have $\hat{x}^*_{hdpr}=0, \forall p\in\mc{P}$ due to \eqref{eq: sors_x-y}. This in turn yields all $z^*_p = 0$ in the optimal subproblem solution due to \eqref{eq: sors_sub_z-x} and thus $\bar{Q}^{s*}_{hdr}=0$. Replacing those values in the cut, we get the following:
\begin{align}
&\overbrace{\bar{Q}^{s*}_{hdr}}^{= 0} \geq \bar{Q}^s_{hdr} - \overbrace{\sum_{p \in \hat{\mc{P}}_{hdr}} c_p^{\text{cancel}}(1 - x^*_{hdpr})}^{= \sum_{p \in \hat{\mc{P}}_{hdr}} c_p^{\text{cancel}} }\label{eq: proof1}\tag{$\star$}
\end{align}

As $\bar{Q}^s_{hdr}$ is the cancellation cost for the $(h, d, r, s)$ tuple corresponding to the patient list $\hat{\mc{P}}_{hdr}$, we have 
\[
\bar{Q}^s_{hdr} \leq \sum_{p \in \hat{\mc{P}}_{hdr}} c_p^{\text{cancel}}  
\]
which follows from the fact that the cancellation cost of an OR cannot exceed the total cancellation cost of all patients assigned to it. We can now conclude that the RHS of \eqref{eq: proof1} is nonpositive. Therefore, the global optimal solution $(\hat{u}^*_{hd}, \hat{y}^*_{hdr}, \hat{x}^*_{hdpr}, \hat{w}^*_p, \bar{Q}^{s*}_{hdr})$ satisfies the LBBD optimality cut. 

\underline{\it Case 2}: If $\hat{y}^*_{hdr} = 1$, then we must have some patients assigned to the $(h, d, r)$ tuple, otherwise we can close this OR and save the cost. At the optimal patient assignment, $\hat{x}^*_{hdpr}$, we either still have all the patients in the set $\hat{\mc{P}}_{hdr}$ in the $(h,d,r)$ tuple, or some patients are no longer assigned to this $(h,d,r)$ tuple. We further discuss those two cases separately:

~~\underline{\it Subcase a}: If all the patients in $\hat{\mc{P}}_{hdr}$ are assigned to the current $(h, d, r)$ tuple at an optimal solution, i.e., $\hat{x}^*_{hdpr} = 1, \forall p \in \hat{\mc{P}}_{hdr}$, then the corresponding cancellation cost for the current patient list, $\bar{Q}^{s*}_{hdr}$, should not be lower than the cancellation cost for $\hat{\mc{P}}_{hdr}$, which is $\bar{Q}^{s}_{hdr}$. That is, we have $\bar{Q}^{s*}_{hdr}\geq\bar{Q}^s_{hdr}$. Substitute the optimal solution into the LBBD cut:
\begin{align*}
&\overbrace{\bar{Q}^{s*}_{hdr}}^{\geq \bar{Q}^s_{hdr}} \geq \bar{Q}^s_{hdr} - \sum_{p \in \hat{\mc{P}}_{hdr}} c_p^{\text{cancel}}\overbrace{(1 - \hat{x}^*_{hdpr})}^{= 0}
\end{align*}
which holds. 

~~\underline{\it Subcase b}: If at the global optimal solution only a subset of patients in $\hat{\mc{P}}_{hdr}$ are still assigned to the current $(h,d,r)$, the proof is more involved. For the ease of proof, we introduce some more notations. Let us denote the patients from $\hat{\mc{P}}_{hdr}$ who are still assigned by set $\hat{\mc{P}}^{A*}_{hdr} \subset \hat{\mc{P}}_{hdr}$, and patients in $\hat{\mc{P}}_{hdr}$ who are no longer assigned as $\hat{\mc{P}}^{N*}_{hdr} =  \hat{\mc{P}}_{hdr}\sm \hat{\mc{P}}^{A*}_{hdr}$. Also, in the global optimal solution there may exist patients who are assigned to the current $(h,d,r)$ but do not belong to $\hat{\mc{P}}_{hdr}$, we denote those patients by $\tilde{\mc{P}}^{A*}_{hdr} \subseteq \mc{P}_{hdr} \sm \hat{\mc{P}}_{hdr}$. Then the set of assigned patients in the global optimal solution is $\mc{P}^{A*}_{hdr} = \hat{\mc{P}}^{A*}_{hdr} \cup \tilde{\mc{P}}^{A*}_{hdr}$. As noted before, the optimal cancellation cost corresponding to the assignment of $\mc{P}^{A*}_{hdr}$ is $\bar{Q}^{s*}_{hdr}$, while the cancellation cost corresponding to $\hat{\mc{P}}_{hdr}$ is $\bar{Q}^s_{hdr}$. It is also useful to find the cancellation cost when only patients in $\hat{\mc{P}}^{A*}_{hdr}$ are assigned to the $(h,d,r,s)$ tuple, which we denote by $\bar{Q}^{s}_{hdr}(\hat{\mc{P}}^{A*}_{hdr})$. It is easy to see that $\bar{Q}^{s*}_{hdr} \geq \bar{Q}^{s}_{hdr}(\hat{\mc{P}}^{A*}_{hdr})$, as $\hat{\mc{P}}^{A*}_{hdr}$ is a subset of $\mc{P}^{A*}_{hdr}$. 

To help illustrate the relationships between patient sets, we use the following simple example in Figure \ref{fig: illuestrate_patient_sets}. There is a set of six patients $\{p_1, ..., p_6\} \in \mc{P}$. The first three of those patients are scheduled by the current master solution to the current $(h,d,r)$, i.e. $\{p_1, p_2, p_3\}\in\hat{\mc{P}}_{hdr}$. In the figure they are  colored with gray. The global optimal solution schedules the patient set $\mc{P}^{A*}_{hdr} = \{p_2, p_3, p_4,p_5\}$.  Then according to our definition, $\hat{\mc{P}}^{A*}_{hdr} = \mc{P}^{A*}_{hdr}\cap\hat{\mc{P}}_{hdr} = \{p_2, p_3\}$, $\tilde{\mc{P}}^{A*}_{hdr} = \mc{P}^{A*}_{hdr}\sm \hat{\mc{P}}^{A*}_{hdr} = \{p_4, p_5\}$, and $\hat{\mc{P}}^{N*}_{hdr} = \hat{\mc{P}}_{hdr}\sm\hat{\mc{P}}^{A*}_{hdr} = \{p_1\}$.

\begin{figure}[h]
\centering
\vspace{0.1cm}
\begin{tikzpicture}[main_node/.style={circle, fill=white!80,draw,inner sep=0pt, minimum size=16pt},	line width=1.2pt]
\tikzset{
ncbar angle/.initial=90,
ncbar/.style={
to path=(\tikztostart)
-- ($(\tikztostart)!#1!\pgfkeysvalueof{/tikz/ncbar angle}:(\tikztotarget)$)
-- ($(\tikztotarget)!($(\tikztostart)!#1!\pgfkeysvalueof{/tikz/ncbar angle}:(\tikztotarget)$)!\pgfkeysvalueof{/tikz/ncbar angle}:(\tikztostart)$)
-- (\tikztotarget)
},
ncbar/.default=-0.5cm,
}
\tikzset{mybrace/.style = {decorate,decoration={brace,amplitude=9pt}}}

\node[main_node, fill = gray!30] (p1) at (0.5,0.8) {$p_1$};
\node[main_node, fill = gray!30] (p2) at (1.5,0.8) {$p_2$};
\node[main_node, fill = gray!30] (p3) at (2.5,0.8) {$p_3$};
\node[main_node] (p4) at (3.5,0.8) {$p_4$};
\node[main_node] (p5) at (4.5, 0.8) {$p_5$};
\node[main_node] (p6) at (5.5,0.8) {$p_6$};
\draw[mybrace, decoration={mirror, amplitude=7pt}] (0.1,0.6) to (0.9,.6) node [midway,xshift=0.5cm] {\small $\hat{\mc{P}}^{N*}_{hdr}$};
\draw[mybrace, decoration={mirror, amplitude=9pt}] (1.1,0.6) to (2.9,.6) node [midway,xshift=2cm] {\small $\hat{\mc{P}}^{A*}_{hdr}$};
\draw[mybrace, decoration={mirror, amplitude=9pt}] (3.1,0.6) to (4.9,.6) node [midway,xshift=4cm] {\small $\tilde{\mc{P}}^{A*}_{hdr}$};
\draw[mybrace, decoration={amplitude=9pt}, yshift = 0.4cm] (1.1,0.6) to (4.9,.6) node [midway,xshift=3cm, yshift = 1.2cm] {\small $\mc{P}^{A*}_{hdr}$};
\end{tikzpicture}

\caption{Illustration of relationships between patient sets}
\label{fig: illuestrate_patient_sets}
\end{figure}
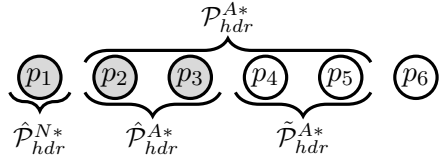

We claim that $\bar{Q}^{s}_{hdr}(\hat{\mc{P}}^{A*}_{hdr}) \geq \bar{Q}_{hdr}^{s} - \sum_{p\in\hat{\mc{P}}^{N*}_{hdr}} c^{\text{cancel}}_p$. To prove this, assume towards contradiction that it is not true. Then we have $\bar{Q}^{s}_{hdr}(\hat{\mc{P}}^{A*}_{hdr}) < \bar{Q}_{hdr}^{s} - \sum_{p\in\hat{\mc{P}}^{N*}_{hdr}} c^{\text{cancel}}_p$, meaning that if in the current OR we have all patients from $\mc{P}^{A*}_{hdr}$, then some patients in $\mc{P}^{N*}_{hdr}$ are also scheduled to this OR, the cancellation cost can increase for at most $\sum_{p\in\hat{\mc{P}}^{N*}_{hdr}} c^{\text{cancel}}_p$. If that is true, then the cancellation cost for the assignment $\hat{\mc{P}}^{A*}_{hdr} \cup \hat{\mc{P}}^{N*}_{hdr}$ is at most $\bar{Q}^{s}_{hdr}(\hat{\mc{P}}^{A*}_{hdr}) + \sum_{p\in\hat{\mc{P}}^{N*}_{hdr}} c^{\text{cancel}}_p < \bar{Q}_{hdr}^{s}$. However, the patient set $\hat{\mc{P}}^{A*}_{hdr} \cup \hat{\mc{P}}^{N*}_{hdr}$ is equivalent to the assignment with patients in the set $\hat{P}_{hdr}$, and its corresponding cancellation cost is exactly $\bar{Q}_{hdr}^{s}$. This is a contradiction.

As the lower bound of $\bar{Q}^{s}_{hdr}(\hat{\mc{P}}_{hdr}^{A*})$ is $\bar{Q}_{hdr}^{s} - \sum_{p\in\hat{\mc{P}}^{N*}_{hdr}} c^{\text{cancel}}_p$, we have the following evaluation of the LBBD cut at the global optimal optimization:
\begin{align*}
&\bar{Q}^{s*}_{hdr} \geq \bar{Q}^{s}_{hdr}(\hat{\mc{P}}_{hdr}^{A*}) \geq \bar{Q}^s_{hdr} - \sum_{p\in\hat{\mc{P}}^N_{hdr}}c^{\text{cancel}}_p\overbrace{(1-\hat{x}^*_{hdp})}^{= 1} - \sum_{p\in \hat{\mc{P}}^A_{hdr}}c^{\text{cancel}}_p \overbrace{(1-\hat{x}^*_{hdp})}^{= 0}
\end{align*}
which is satisfied thanks to the relations of the cut's LHS and RHS to the middle comparative term as mentioned above.\qed
\endproof

\subsection{Proof of Theorem \ref{theorem: qlb}} \label{subch: app_proof_th2}
{\cred\begin{theorem}
Constraints \eqref{eq: qlb}, which are defined as 
\begin{align}
Q^s_{hdr} \geq \Bigg(\min_{p\in\mc{P}} \frac{c_p^{\text{cancel}}}{T^s_p}\Bigg) \Bigg(\sum_{p\in\mc{P}} T^s_p x_{hdpr} - B_{hd}\Bigg) ~~\forall h\in\mc{H}, d\in\mc{D}, r\in\mc{R}_h, s\in\mc{S},\nonumber
\end{align}
\end{theorem}
provide valid lower bounds (LBs) for $Q^s_{hdr}$.}
\proof{Proof.}
At optimality $Q^s_{hdr}$ equals the optimal objective value of the subproblem \eqref{model: sors_sub}. We want to prove that the RHS of inequality \eqref{eq: qlb} is either trivially true or otherwise can be obtained by relaxing the subproblem. 

First, we look at the trivial case where the OR corresponding to $Q^s_{hdr}$ is not opened. In this case, there should be no cost for cancelling as no patient is assigned in the first place. Therefore, $Q^s_{hdr} = 0$. Since in this case the RHS of the inequality becomes $\min_{p\in\mc{P}} \bigg(\frac{c_p^{\text{cancel}}}{T^s_p}\bigg)(- B_{hd}) < 0$ as the consequence of $x_{hdpr} = 0~(\forall p\in\mc{P})$, the constraint is valid.

If the OR corresponding to $Q^s_{hdr}$ is opened in an optimal solution, there must exist at least one patient who is assigned to this OR. Given an assignment of patients, $\hat{x}_{hdpr}~(\forall p\in\mc{P})$, and the set of assigned patients, $\hat{\mc{P}}_{hdr}$, suppose a subset of the patients, $\hat{\mc{P}}^C_{hdr}\subseteq\hat{\mc{P}}_{hdr}$, is cancelled as dictated by the optimal solution of the subproblem \eqref{model: sors_sub}. Then by definition we have $Q^s_{hdr} = \sum_{p \in \hat{\mc{P}}^C_{hdr}} c_{p}^{\text{cancel}}$. Due to the fact that after cancellation, the total surgery duration of accepted patients should be no more than the operating time limit, $B_{hd}$, we have the following:
\begin{align*}
\sum_{p \in \hat{\mc{P}}^C_{hdr}} c_{p}^{\text{cancel}} & = \sum_{p \in \hat{\mc{P}}^C_{hdr}} c_{p}^{\text{cancel}} \hat{x}_{hdpr}\\
& = \sum_{p \in \hat{\mc{P}}^C_{hdr}} \frac{c_{p}^{\text{cancel}}}{T^s_p} T^s_p \hat{x}_{hdpr}\\
&\geq\left(\min_{p\in \mc{P}} \frac{c_{p}^{\text{cancel}}}{T^s_p}\right) \sum_{p \in \hat{\mc{P}}^C_{hdr}} T^s_p \hat{x}_{hdpr}\\
&\geq \left(\min_{p\in \mc{P}} \frac{c_{p}^{\text{cancel}}}{T^s_p}\right) \left(\sum_{p \in \hat{\mc{P}}^C_{hdr}} T^s_p \hat{x}_{hdpr} + \left(\sum_{p\in \mc{P}\sm\hat{\mc{P}}^C_{hdr}} T^s_p \hat{x}_{hdpr} - B_{hd}\right)\right)\\
& = \left(\min_{p\in \mc{P}} \frac{c_{p}^{\text{cancel}}}{T^s_p}\right)\left(\sum_{p\in\mc{P}} T^s_p \hat{x}_{hdpr} - B_{hd}\right)
\end{align*}

Therefore, for any assignment $\hat{x}_{hdpr}~(\forall p\in\mc{P})$, the expression $\left(\min_{p\in \mc{P}} \frac{c_{p}^{\text{cancel}}}{T^s_p}\right)\left(\sum_{p\in\mc{P}} T^s_p \hat{x}_{hdpr} - B_{hd}\right)$ provides a LB for $Q^s_{hdr}$. Substitute the fixed assignment $\hat{x}_{hdpr}$ with the variable $x_{hdpr}$ and we get the constraint \eqref{eq: qlb}.\qed
\endproof
\subsection{Proof of Theorem \ref{theorem: lbbd2_valid}}\label{subch: app_proof_th3}
{\cred\begin{theorem}
The LBBD optimality cut \eqref{eq: lbbd_cut2}, which is defined as 
 \begin{subequations}
\begin{alignat}{5}
&Q_{hd}\geq \bar{Q}_{hd}\left(g_{hdj}-\sum_{p\in\hat{\mc{P}}_{hd}}(1-x_{hdp})\right)\nonumber\\
&y_{hd}\geq (1+\hat{y}_{hd})(1-g_{hdj}),\nonumber\\
&g_{hdj}\in\{0,1\}\nonumber
\end{alignat}
\end{subequations}
is valid, where $\bar{Q}_{hd}$ is the optimal objective value of the LBBD subproblem \eqref{model: sors2_sub}, and $\hat{y}_{hd}$ is the optimal solution of $y_{hd}$ from the LBBD master problem \eqref{model: sors2_lbbd_master}.
\end{theorem}}
\proof{Proof.}
For any $(h,d)$ pair, we need to prove that any global optimal solution $(\hat{u}'_{hd}, \hat{y}'_{hd}, \hat{x}'_{hdp}, \bar{Q}'_{hd})$ is not excluded by the optimality cut. Given an optimal \DE~solution, $(\hat{u}'_{hd}, \hat{y}'_{hd}, \hat{x}'_{hdp})$ denotes the corresponding solutions for the main LBBD master decisions, and $\bar{Q}'_{hd}$ is the corresponding optimal objective of the subproblem. \\
\underline{\it Case 1}: If $\hat{u}'_{hd}=0$, then the only global optimal solution for the current $h$ and $d$ will be $\hat{y}'_{hd}=0, \hat{x}'_{hdp}=0, \bar{Q}'_{hd}=0$. This is feasible to the optimality cut. \\
\underline{\it Case 2}: If $\hat{u}'_{hd}=1$, then we discuss the following two cases separately:\\
\indent \underline{\it Subcase a}: If $\hat{y}'_{hd}>\hat{y}_{hd}$, we are not able to find a nontrivial LB for $Q_{hd}$ without solving another subproblem, because when there are more ORs available, the cancellation cost of the current $(h,d)$ pair can be either zero or some nonzero value that is smaller than $\bar{Q}_{hd}$. This is why we make the cut redundant in this case: let $g_{hdj}=0$, \eqref{eq: LBBD3_1} becomes $$Q_{hd}\geq\bar{Q}_{hd}\left(-\sum_{p\in\hat{P}_{hd}}(1-x_{hdp})\right)$$
This is always true since RHS is nonpositive, so  $\bar{Q}'_{hd}$ and $\hat{x}'_{hdp}$ also satisfy this inequality. \eqref{eq: LBBD3_2} is now $y_{hd}\geq(1+\hat{y}_{hd})$, which is equivalent to $\hat{y}_{hd}>\hat{y}_{hd}$, and that is exactly the assumption of this case. \\
\indent \underline{\it Subcase b}: If $\hat{y}'_{hd}\leq\hat{y}_{hd}$, then there are further two subcases to discuss. Note that due to \eqref{eq: LBBD3_2}, we always have $g_{hdj} = 1$ in this case.\\
\indent~\underline{\it Subcase b1}: If in the global optimal all patients that are assigned in the current solution are still assigned, i.e. $\hat{x}'_{hdp}\geq\hat{x}_{hdp}, \forall p\in \mc{P}_{hd}$, then we claim that the optimal cancellation cost will not decrease from the current value $\bar{Q}_{hd}$ because there are the same number of or a smaller number of rooms, but all the current assigned patients are still assigned. To see this argument is another way, assume for contradiction that $Q_{hd} < \bar{Q}_{hd}$ in this case. Then this means we can also schedule the patients in the current patient list in the same way with a lower cancellation cost than $\bar{Q}_{hd}$, which is a contradiction. Therefore, $Q_{hd} \geq \bar{Q}_{hd}$. Let $g_{hdj}=1$, and also replace $Q_{hd}$ and $x_{hdp}$ with $\hat{Q}'_{hdp}$ and $\hat{x}'_{hdp}$, \eqref{eq: lbbd_cut2} becomes: 
\begin{align*}
&\overbrace{\hat{Q}'_{hd}}^{\geq \bar{Q}_{hd}} \geq \bar{Q}_{hd}\left(\overbrace{g_{hdj}}^{=1}-\overbrace{\sum_{p\in\hat{P}_{hd}}(1-\hat{x}'_{hdp})}^{=0}\right)\\
&\overbrace{\hat{y}'_{hd}}^{\leq \hat{y}_{hd}}\geq (1+\hat{y}_{hd})(1-\overbrace{g_{hdj}}^{=1}) 
\end{align*}
Therefore, in this case the global solution also satisfies the optimality cut.\\
\indent~\underline{\it Subcase b2}: If some currently assigned patients are no longer assigned, then we cannot give a nontrivial LB for the global optimal $\bar{Q}'_{hd}$, because the cancellation cost can either be zero or be a nonzero value that is larger or smaller than $\bar{Q}_{hd}$. Thus we make the LBBD cut redundant:
\begin{align*}
&\overbrace{\hat{Q}'_{hd}}^{\geq 0} \geq \bar{Q}_{hd}\left(\overbrace{g_{hdj}}^{=1}-\overbrace{\sum_{p\in\hat{P}_{hd}}(1-\hat{x}'_{hdp})}^{\geq 1}\right)\\
&\overbrace{\hat{y}'_{hd}}^{\leq \hat{y}_{hd}}\geq (1+\hat{y}_{hd})(1-\overbrace{g_{hdj}}^{=1}) 
\end{align*}
which is satisfied by the global optimal solution. \qed
\endproof
\subsection{Proof of Theorem \ref{theorem: qlb2}}\label{subch: app_proof_th4}
{\cred\begin{theorem}
The constraints \eqref{eq: qlb2}, which are defined as
\begin{align}
Q_{hd} \geq \frac{1}{|\mc{S}|} \sum_{s\in\mc{S}} \left(\min_{p\in\mc{P}} \frac{c_p^{\text{cancel}}}{T^s_p}\right)\left(\sum_{p\in\mc{P}} T^s_p x_{hdp} - B_{hd} y_{hd}\right) ~~~\forall h\in\mc{H}, d\in\mc{D},\nonumber
\end{align}
\end{theorem}
are valid for (DE).}
\proof{Proof.}
At optimality $Q_{hd}$ equals the optimal objective value of the subproblem \eqref{model: sors2_sub}. We want to prove that the RHS of inequality \eqref{eq: qlb2} is either trivially true or otherwise can be obtained by relaxing the subproblem. 

First, we look at the trivia case where the $(h,d)$ pair corresponding to $Q_{hd}$ is not opened. In this case, there should be no cost for cancelling as no patient is assigned in the first place. Therefore, $Q_{hd} = 0$. Since in this case the RHS of the inequality becomes 0 as the consequence of $x_{hdp} = 0~(\forall p\in\mc{P})$ and $y_{hd} = 0$, the constraint is valid.

If the $(h,d)$ pair corresponding to $Q_{hd}$ is opened in an optimal solution, there must exist at least one patient who is assigned to this $(h,d)$ pair. Given an assignment of patients, $\hat{x}_{hdp}~(\forall p\in\mc{P})$, the set of assigned patients, $\hat{\mc{P}}_{hd}$, and the number of opened ORs, $\hat{y}_{hd}$. Suppose under the scenario $s$ a set of patients, $\hat{\mc{P}}^{sC}_{hdr}\subseteq\hat{\mc{P}}_{hd}$, is cancelled as dictated by the optimal solution of the subproblem \eqref{model: sors2_sub}. Then by definition we have $Q_{hd} =\frac{1}{|\mc{S}|}\sum_{s\in\mc{S}} \sum_{p \in \hat{\mc{P}}^{sC}_{hdr}} c_{p}^{\text{cancel}}$. Due to the fact that after cancellation, the total surgery duration of accepted patients should be no more than the total operating time limit of opened ORs, $B_{hd}\hat{y}_{hd}$, we have the following for any scenario $s\in\mc{\mc{S}}$:
\begin{align*}
\sum_{p \in \hat{\mc{P}}^{sC}_{hdr}} c_{p}^{\text{cancel}} & = \sum_{p \in \hat{\mc{P}}^{sC}_{hdr}} c_{p}^{\text{cancel}} \hat{x}_{hdp}\\
& = \sum_{p \in \hat{\mc{P}}^{sC}_{hdr}} \frac{c_{p}^{\text{cancel}}}{T^s_p} T^s_p \hat{x}_{hdp}\\
&\geq\left(\min_{p\in \mc{P}} \frac{c_{p}^{\text{cancel}}}{T^s_p}\right) \sum_{p \in \hat{\mc{P}}^{sC}_{hdr}} T^s_p \hat{x}_{hdp}\\
&\geq \left(\min_{p\in \mc{P}} \frac{c_{p}^{\text{cancel}}}{T^s_p}\right) \left(\sum_{p \in \hat{\mc{P}}^{sC}_{hdr}} T^s_p \hat{x}_{hdp} + \left(\sum_{p\in \mc{P}\sm\hat{\mc{P}}^{sC}_{hdr}} T^s_p \hat{x}_{hdp} - B_{hd}\hat{y}_{hd}\right)\right)\\
& = \left(\min_{p\in \mc{P}} \frac{c_{p}^{\text{cancel}}}{T^s_p}\right)\left(\sum_{p\in\mc{P}} T^s_p \hat{x}_{hdp} - B_{hd}\hat{y}_{hd}\right)
\end{align*}

Therefore, for any assignment $\hat{x}_{hdp}~(\forall p\in\mc{P})$, the expression $\frac{1}{|\mc{S}|} \sum_{s\in\mc{S}} \bigg(\min_{p\in\mc{P}} \frac{c_p^{\text{cancel}}}{T^s_p}\bigg)\bigg(\sum_{p\in\mc{P}} T^s_p \hat{x}_{hdp} - B_{hd} \hat{y}_{hd}\bigg)$ provides a LB for $Q_{hd}$. Substitute the fixed assignment $\hat{x}_{hdp}$ with the variable $x_{hdp}$ and $\hat{y}_{hd}$ with $y_{hd}$ then we get the constraint \eqref{eq: qlb2}.
\qed
\endproof

\section{Flow Chart for the Two-stage Decomposition}\label{chapter: app_flowchart}
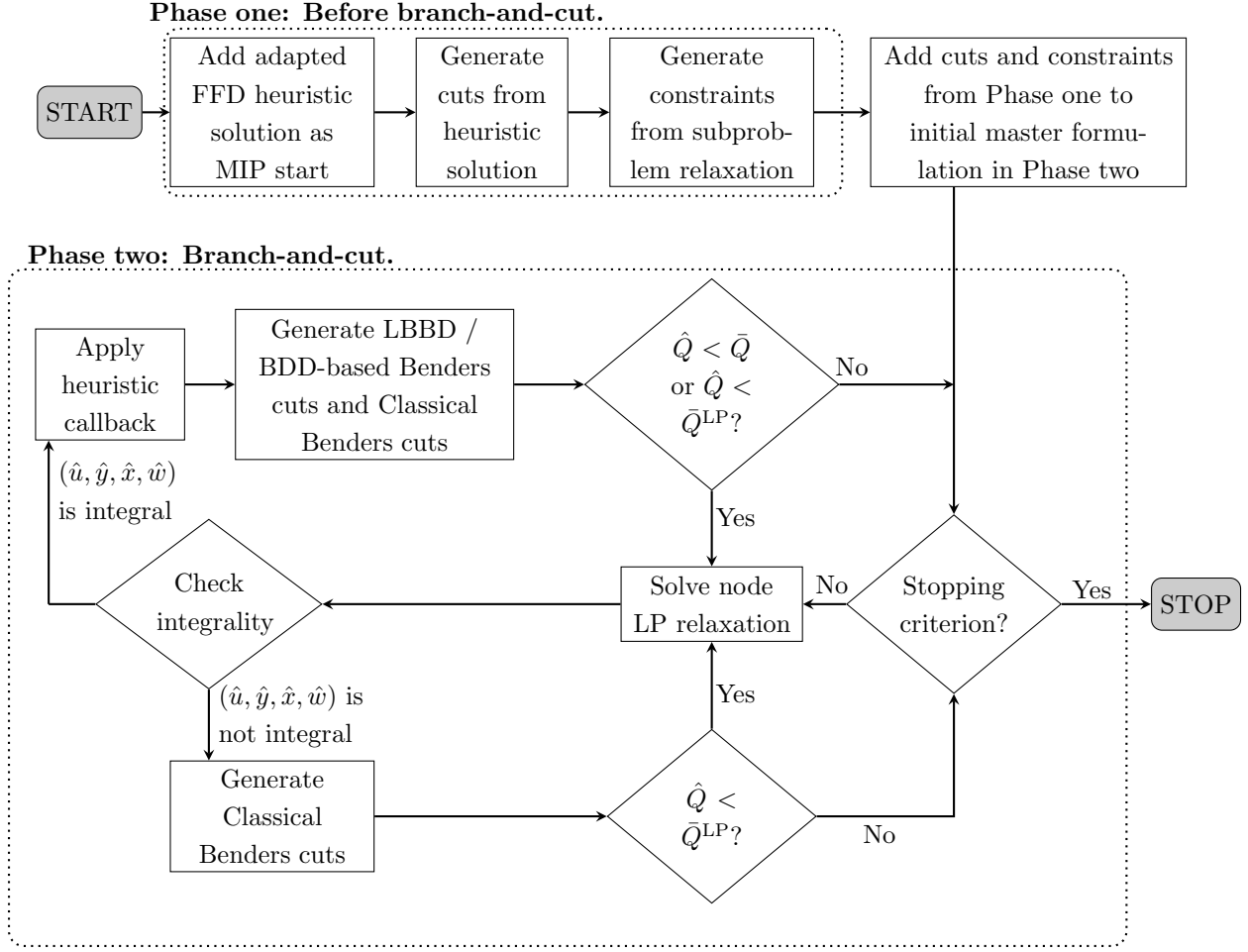
\begin{figure}[H] 
\begin{center}
{
\begin{tikzpicture}[node distance=1.75cm]
\tikzstyle{startstop} = [rectangle, rounded corners, minimum width=1.15cm, minimum height=0.7cm,text centered, draw=black, fill = gray!40!white]
\tikzstyle{process} = [rectangle, minimum width=2cm, minimum height=1cm, text centered, text width=1.98cm, draw=black]
\tikzstyle{decision} = [diamond, aspect=1.4, text centered, text width=1.2cm, draw=black]
\tikzstyle{arrow} = [thick,->,>=stealth]
\small
\node (start) [startstop] {START};
\node (lp1) [process, right of=start,xshift = 0.7cm,text width=2.5cm] {Add adapted FFD heuristic solution as MIP start};
\node (ben1) [process, right of=lp1,xshift = 1.2cm,text width=1.8cm] {Generate cuts from heuristic solution};
\node (vio1) [process, right of=ben1,xshift = 1.2cm,text width=2.5cm] {Generate\\ constraints from subproblem relaxation};

\node (trans) [process, right of=vio1,text width=4cm,xshift = 2.5cm] {Add cuts and constraints \\ from Phase one to initial master formulation in Phase two};
\node (dummytrans) [xshift = -1cm] at (trans.south) {};

\node (stopping) [decision, below of=dummytrans,aspect=1.2,yshift = -3.85cm, text width=1.47cm] {Stopping criterion?};
\node (stop) [startstop, right of=stopping,xshift = 1.5cm] {STOP};
\node (lp2) [process,text width=2.2cm] at (vio1 |- stopping) {Solve node LP relaxation};

\node (vio2) [decision, above of=lp2,aspect=1.2,yshift = 1.2cm] {$\hat{Q} < \bar{Q}$ or $\hat{Q} < \bar{Q}^{\text{LP}}$?};
\node (vio3) [decision, below of=lp2,aspect=1.2,yshift = -1.1cm] {$\hat{Q} < \bar{Q}^{\text{LP}}$?};


\node (ben2) [process,text width=3.5cm] at (lp1.east |- vio2) {Generate LBBD / BDD-based Benders cuts and Classical \\Benders cuts};
\node (ben3) [process,text width=2.5cm] at (lp1 |- vio3) {Generate Classical Benders cuts};

\node (UBheur) [process,left of=ben2,xshift= -1.8cm,text width=1.5cm, text centered] {Apply heuristic callback};

\node (check) [decision, aspect=1.33,text width=1.4cm ,left of=lp2,xshift= -5cm] {Check \\ integrality};
\node (dummycheck1) [left of=check, xshift=-0.4cm] {};
\node (dummyinf) [left of=check, xshift=-0.48cm] {};
\node (dummycheck2) [below of=dummyinf, yshift=-2.6cm] {};
\node (dummy3) [right of=vio3, xshift=0.35cm] {};

\draw [arrow] (start) -- (lp1);
\draw [arrow] (lp1) -- (ben1);
\draw [arrow] (ben1) -- (vio1);
\draw [arrow] (vio1) -- (trans);
\draw [arrow] (dummytrans.center) -- (stopping);
\draw [arrow] (stopping) -- node[anchor=east,xshift = 0.2cm,yshift=0.2cm] {Yes} (stop);
\draw [arrow] (stopping) -- node[anchor=south,xshift = 0.1cm,yshift=0.01cm] {No} (lp2);

\draw [arrow] (vio2) -- node[anchor=south,xshift = 0.32cm,yshift=-0.1cm] {Yes} (lp2);
\draw [arrow] (vio3) -- node[anchor=south,xshift = 0.32cm,yshift=-0.36cm] {Yes} (lp2);


\draw [arrow] (vio2) -- node[anchor=south,xshift = -0.6cm,yshift=-0.01cm] {No} (vio2 -| stopping);
\draw [arrow] (vio3) -| node[anchor=north,xshift = -1 cm,yshift=0.01cm] {No} (stopping);

\draw [arrow] (ben2) -- (vio2);
\draw [arrow] (ben3) -- (vio3);

\draw [arrow] (UBheur) -- (ben2);

\draw [arrow] (lp2) -- (check);
\draw [arrow] (dummycheck1.center) -- node[anchor=west,text width=1.85cm,xshift = 0cm,yshift=0.4cm] {$(\hat{u}, \hat{y}, \hat{x}, \hat{w})$ is integral } (UBheur.south -| dummycheck1.center);
\draw [arrow] (check.south) -- node[anchor=west,text width=2.3cm,xshift = 0cm,yshift=0.12cm] {$(\hat{u}, \hat{y}, \hat{x}, \hat{w})$ is not integral} (check.south |- ben3.north);

\draw [-,thick] (check) -- (dummycheck1.center);

\node [thick,xshift=0.15cm,draw=black,rounded corners=2mm, dotted, minimum width=9.1cm, minimum height=0.1cm, fit = (lp1) (vio1) ] (Phase1)  {};
\node [thick,xshift=0.27cm,draw=black,rounded corners=2mm, dotted, minimum width=15cm, minimum height=8.2cm, fit = (UBheur) (dummycheck2) (stopping) (vio2.north)] (Phase2)  {};

\node (p1) [above of=Phase1,xshift = 0.25cm, yshift = -0.4cm,text width=10cm] {\textbf{Phase one: Before branch-and-cut.}};
\node (p2) [above of=Phase2,xshift = -4.75cm, yshift = 3cm,text width=5cm] {\textbf{Phase two: Branch-and-cut.}};
\end{tikzpicture}
}
\caption{Flow chart of the two-stage decomposition algorithm.}
\label{fig: impl_flowchart1}
\end{center}
\end{figure}
\section{Classical Benders Cuts for the Three-stage Decomposition}\label{chapter: app_benders3}
For the LBBD subproblem, we relax the LBBD subproblem \eqref{model: sors2_sub}. The binary variables $x_{pr}$, $y_{r}$, and $z_{pr}^s$ are redefined as continuous variables. The variables in the parentheses at the end of constraints some constraints are their corresponding dual variables. We denote the {\cred relaxed LBBD subproblem by $\mc{Q}_{hd}^{\text{LP}}(\hat{x}_{hd\cdot}, \hat{y}_{hd}, T^s_\cdot)$, and shorten it as $\bar{Q}_{hd}^{\text{LP}}$ later in the text}:

\begin{align*}
{\cred \mc{Q}_{hd}^{\text{LP}}(\hat{x}_{hd\cdot}, \hat{y}_{hd}, T^s_\cdot)} = \min~~& \frac{1}{|\mc{S}|}\sum_{s\in\mc{S}}\sum_{p\in \mc{P}}\sum_{r\in\mc{R}_{h}} c_p(x_{pr}-z^s_{pr}) \hspace*{0cm}\\
\st ~~& \sum_{r\in\mc{R}_h}x_{pr}=\hat{x}_{hdp} & & \forall p \in \mc{P} &&&(\gamma_{p})\\
     & \sum_{p\in \mc{P}}T^s_p z^s_{pr} \leq B_{hd}y_r  & & \forall r\in\mc{R}_h, s\in\mc{S} \\
     & z^s_{pr} \leq x_{pr} & & \forall p\in \mc{P}, r\in\mc{R}_h, s\in\mc{S}\\
     & x_{pr} \leq y_{r} & &\forall p\in \mc{P}, r\in\mc{R}_h\\
     & \sum_{r\in\mc{R}_h} y_{r} \leq \hat{y}_{hd} && &&&(\beta) \\
     &y_r\leq 1& &\forall r\in\mc{R}_h &&~&(\delta_{r}) \\
     & y_r, x_{pr},z^s_{pr}\geq 0& &\forall p\in \mc{P}, r\in\mc{R}_h, s\in\mc{S}
\end{align*}

When $\hat{Q}_{hd} <\bar{Q}^{\text{LP}}_{hd}$, the following classical Benders cuts are added to the LBBD master problem:
\begin{align*}\label{eq: lbbd_bd_cut2}
Q_{hd}\geq \sum_{p\in \mc{P}}\bar{\gamma}_{p}x_{hdp}+\bar{\beta}y_{hd}+\sum_{r\in\mc{R}_h}\bar{\delta}_{r}
\end{align*}
where $\bar{\gamma}_p$, $\bar{\beta}$, and $\bar{\delta}_r$ are the optimal solutions of their corresponding dual variables. 

For the decomposition of LBBD subproblem, we relax the variable $z^s_{pr}$ as a continuous variable in subproblem \eqref{model: bdd2_sub}. {\cred We denote the relaxed subproblem by $\theta_{sr}^{\text{LP}}(\check{x}_{\cdot r}, T^s_\cdot)$, and shorten it as $\ddot{\theta}_{sr}^{\text{LP}}$ later in the text}:
\begin{align*}
{\cred \theta_{sr}^{\text{LP}}(\check{x}_{\cdot r}, T^s_\cdot)}= \min~~& \sum_{p\in \mc{P}} - c_p z^s_{pr} \hspace*{0cm}\\
\st ~~& \sum_{p\in \mc{P}} T^s_p z^s_{pr} \leq B_{hd} && &&&(\eta)  \\
     & z^s_{pr} \leq \check{x}_{pr} & & \forall p\in \mc{P}~~~ &&&(\iota)\\
     & z^s_{pr}\in \{0,1\}& &\forall p\in \mc{P}
\end{align*}

The corresponding classical Benders cut is:
\begin{align*}
\theta_{sr} \geq \sum_{p\in\mc{P}}\ddot{\i_p} x_{pr}+B_{hd}\ddot{\eta}
\end{align*}
where $\ddot{\eta}$ and $\ddot{\i}_p$ are the optimal values for the corresponding dual variables.
\section{Overall Implementation Approach for Three-stage Decomposition}\label{chapter: app_impl3}
In this section we first describe the LBBD decomposition, then explain the decomposition of LBBD subproblem. 

{\bf LBBD decomposition: }

\underline{\it Phase one}: This phase is very similar to the phase one of two-stage decomposition in Section \ref{subch: overall_implement}. We first use the adapted FFD heuristic to obtain an initial solution. This solution is added as a warm start in the commercial solver to provide a feasible solution at the start of branch-and-cut. We also generate LBBD cuts and classical Benders cuts from this solution and add them to the LBBD master problem. Next, we generated the constraints \eqref{eq: qlb2} from the subproblem relaxations and also add them to the LBBD master problem.  

\underline{\it Phase two}: We obtain the LBBD master problem with extra cuts and constraints from phase one and solve it with branch-and-cut. At each branch-and-bound node solve the node LP relaxation. If the objective value is greater or equal to the incumbent UB, then the current node can be pruned. Otherwise if the objective value is less than the incumbent UB, we proceed to check the integrality of $(\hat{u}, \hat{y}, \hat{x}, \hat{w})$ in the master solution. If $(\hat{u}, \hat{y}, \hat{x}, \hat{w})$ is integral, then solve the corresponding LBBD subproblems \eqref{model: sors2_sub} with further decomposition (described with detail in the next paragraph). In the process of solving the LBBD subproblem, check if we need early stopping as explained in Section \ref{chapter: early_stop}. If the solving process early stops, add LBBD cuts \eqref{eq: ES_LBBD_cuts}; otherwise, solve the LBBD subproblem to get the optimal objective value $\bar{Q}$. Also, solve subproblem LP relaxations and get optimal objective value $\bar{Q}^{\text{LP}}$. Decide if we can insert an additional heuristic solution as described in Section \ref{chapter: heurcb}. Also, generate LBBD cuts and classical Benders cuts. In the CPLEX lazy constraint callback, compare the master solution of $\hat{Q}$ with $\bar{Q}$ and $\bar{Q}^{\text{LP}}$. If $\hat{Q} < \bar{Q}$ then add LBBD cuts; if $\hat{Q} < \bar{Q}^{\text{LP}}$ add classical Benders cuts. On the other hand, if some elements in the solution $(\hat{u}, \hat{y}, \hat{x}, \hat{w})$ are fractional, we only solve the subproblem LP relaxations, obtain $\bar{Q}^{\text{LP}}$, generate the classical Benders cuts and implement them if $\hat{Q} < \bar{Q}$ within the CPLEX user cut callback. We use the same user cut management as in Section \ref{subch: overall_implement} to manage those user cuts. After cutting planes are added in the CPLEX lazy constraint callback or the user cut callback, the node LP relaxation is solved again with those additional cutting planes. We repeat this process, until the stopping criteria is met, i.e. the gap between branch-and-bound UB and LB is small enough. In our implementation we stop the algorithm when such a gap is no more than 1\%.

{\bf Decomposition of LBBD subproblem: }

\underline{\it Phase one}: Use the adapted FFD heuristic to obtain an initial solution. This solution is added as a warm start in the commercial solver to provide a feasible solution at the start of branch-and-cut. 

\underline{\it Phase two}: We solve the BDD master problem with branch-and-cut. At each branch-and-bound node solve the node LP relaxation. If the objective value is greater or equal to the incumbent UB, then the current node can be pruned. Otherwise if the objective value is less than the incumbent UB, we proceed to check the integrality of $\check{x}$ in the master solution. If $\check{x}$ is integral, then solve the corresponding BDD subproblems \eqref{model: bdd2_sub} and the BDD subproblem LP relaxations to get their respective optimal objective values $\ddot{\theta}$ and $\ddot{\theta}^{\text{LP}}$. Generate BDD-based Benders cuts and classical Benders cuts. In the CPLEX lazy constraint callback, compare the master solution of $\check{\theta}$ with $\ddot{\theta}$ and $\ddot{\theta}^{\text{LP}}$. If $\check{\theta} < \ddot{\theta}$ then add BDD-based Benders cuts; if $\check{\theta} < \ddot{\theta}^{\text{LP}}$ add classical Benders cuts. On the other hand, if some elements in the solution $\check{x}$ are fractional, we only solve the subproblem LP relaxations, obtain $\ddot{\theta}^{\text{LP}}$, generate the classical Benders cuts and implement them if $\check{\theta} < \ddot{\theta}$ within the CPLEX user cut callback. The user cut management and stopping criteria are the same as in the LBBD decomposition.
\section{Parameter Values for Computational Analysis}\label{subch: app_table_parameter}
\begin{table}[H]
\small
\centering
\renewcommand\thetable{A.1}\label{table: parameters}
\caption{Parameter values}
\begin{tabular}{ll}
\hline
$\kappa_1$ & 50 dollars                                                                         \\
$\kappa_2$ & -5 dollars                                                                        \\
$\kappa_3$ & -80 dollars                                                                       \\
$\kappa_4$ & -100 dollars                                                                      \\
$\Gamma$   & 500                                                                                \\
$\rho_p$   & Uniform distribution in \{1,2,...,5\}, where 1 is the least urgent 5 is the most urgent \\
$B_{hd}$   & Uniform distribution {[}420, 480{]}  minutes in 15-minute intervals                \\
$\alpha_p$ & Uniform distribution {[}60, 120{]} days                                            \\
$F_{hd}$   & Uniform distribution {[}4000, 6000{]}                                              \\
$G_{hd}$   & Uniform distribution {[}1500, 2500{]}                                             \\
$\cred c^{\text{sched}}_{dp}$ & $\cred \kappa_1\rho_p(d-\alpha_p)$ \\
$\cred c^{\text{unsched}}_{p}$ & $\cred \kappa_2\rho_p(|\mc{D}| + 1 -\alpha_p)$ \\
$\cred c^{\text{cancel}}_{p}$ & $\cred \kappa_3\rho_p(|\mc{D}| + 1 -\alpha_p), \forall p\in\mc{P}\sm\mc{P}^'$ \\
$\cred c^{\text{cancel}}_{p}$ & $\cred \kappa_4\rho_p(|\mc{D}| + 1 -\alpha_p), \forall p\in\mc{P}^'$ \\
$\cred T^s_p$ & \cred Truncated lognormal distribution with mean 160 minutes, standard deviation 40, and truncation \\
	&\cred range [45, 480]\\
$\cred \omega_p$ & $\cred (\alpha_p - |\mc{D}|)\rho_p$\\
\hline
\end{tabular}
\end{table}

\section{Best Integer Solution Results for Algorithm Comparison}\label{ch: app_best_int}
\begin{table}[H]
\centering
\small
\cred
\renewcommand\thetable{A.2} 
\caption{Comparison of Algorithms (continued): Best Integer Solution Objective Values (i.e., Upper Bounds)}
\label{tab:my-table}
\begin{tabular}{crrrr}
\hline
instance  & \multicolumn{1}{l}{}    & \multicolumn{1}{l}{}      & \multicolumn{1}{l}{}       & \multicolumn{1}{l}{}       \\
(p-h-d-r) & \multicolumn{1}{c}{MIP} & \multicolumn{1}{c}{2-BDD} & \multicolumn{1}{c}{2-LBBD} & \multicolumn{1}{c}{3-LBBD} \\ \hline
10-2-3-3  & -117624                 & -117624                   & -117670                    & -117670                    \\
25-2-3-3  & -248582                 & -252010                   & -253491                    & -249238                    \\
10-3-5-3  & -117227                 & -117671                   & -117595                    & -117671                    \\
25-3-5-3  & -241882                 & -247676                   & -247807                    & -248107                    \\
50-3-5-3  & -426982                 & -433702                   & -391874                    & -380486                    \\
75-3-5-3  & -600867                 & -677288                   & -679017                    & -638701                    \\
10-2-3-5  & -119481                 & -119551                   & -118935                    & -118935                    \\
25-2-3-5  & -253762                 & -253335                   & -255321                    & -                          \\
50-2-3-5  & -356952                 & -435996                   & -452404                    & -                          \\
75-2-3-5  & -611189                 & -696665                   & -694526                    & -                          \\
10-3-5-5  & -119555                 & -119537                   & -117611                    & -119588                    \\
25-3-5-5  & -240490                 & -251945                   & -251931                    & -                          \\
50-3-5-5  & -359391                 & -449867                   & -448060                    & -                          \\
75-3-5-5  & -613869                 & -676269                   & -726670                    & -                          \\ \hline
\end{tabular}
\end{table}
{\cred Notice that for instance 10-2-3-5, both 2-BDD and 3-LBBD are solved to optimality, but they have different best integer results. This difference is caused by setting the 1\% relative MIP gap in the solver.}

\end{document}